\documentclass[secnum,leqno]{rims-bessatsu}
\usepackage{amssymb, amsmath, bm, bbm}
\usepackage{graphics}
\usepackage[dvips]{graphicx}
\usepackage{eepic}
\usepackage{epic}
\newtheorem{theorem}{Theorem}[section]
\newtheorem{lemma}[theorem]{Lemma}
\newtheorem{proposition}[theorem]{Proposition}
\newtheorem{corollary}[theorem]{Corollary}

\theoremstyle{definition}
\newtheorem{definition}[theorem]{Definition}
\newtheorem{example}[theorem]{Example}
\newtheorem{remark}[theorem]{Remark}
\newtheorem{question}[theorem]{Question}
\newtheorem*{acknowledgment}{Acknowledgment}
\theoremstyle{remark}
\newtheorem*{notation}{Notation}
\newcommand{\bZ}{\mathbb{Z}}
\newcommand{\bQ}{\mathbb{Q}}
\newcommand{\bC}{\mathbb{C}}\newcommand{\bG}{\mathbb{G}}
\newcommand{\cC}{\mathcal{C}}\newcommand{\cD}{\mathcal{D}}
\newcommand{\cH}{\mathcal{H}}\newcommand{\cS}{\mathcal{S}}
\newcommand{\cT}{\mathcal{T}}

\def\fn#1{\operatorname{#1}} 

\newcommand{\ep}{\varepsilon}

\def\fn#1{\mathop{{\rm #1}\vphantom{\sin}}} 
\def\bm#1{\mathbbm{#1}}

\def\ul#1{\mathop{\underline{#1}}}

\title{Rationality problem for quasi-monomial actions}
\author{Akinari \textsc{Hoshi}\footnote{Department of Mathematics, 
Rikkyo University, 3--34--1 Nishi-Ikebukuro Toshima-ku, Tokyo, 171-8501, Japan.
\newline e-mail: \texttt{hoshi@rikkyo.ac.jp}}}
\AuthorHead{Akinari Hoshi}         
\support{This work was partially supported by KAKENHI (22740028).}

\VolumeNo{51}           
\YearNo{2014}           
\PagesNo{203--227}      
\communication{Received March 25, 2013. 
Revised February 10, 2014, May 21, 2014 and June 1, 2014. 
}      
\begin{document}
%

\maketitle

\begin{abstract}      
We give a short survey of the rationality problem 
for quasi-monomial actions which includes Noether's problem 
and the rationality problem for algebraic tori, 
and report some results on rationality problem 
in three recent papers 
Hoshi, Kang and Kitayama [HKKi], 
Hoshi, Kang and Kunyavskii [HKKu] and 
Hoshi and Yamasaki [HY].
\end{abstract}


\tableofcontents

In Section \ref{se1}, 
we recall first the definitions of rationalities for field extensions, 
concerning which we give a survey in this article. 
In Subsection \ref{sse11}, 
we explain that the rationality problem 
for quasi-monomial actions includes Noether's problem, 
the rationality problem for algebraic tori and 
the rationality problem for Severi-Brauer varieties. 
In Subsection \ref{sse12}, 
we list known results for monomial actions 
in dimension two and three 
due to Hajja, Kang, Saltman, Hoshi, Rikuna, Prokhorov, Kitayama and Yamasaki. 
In Subsection \ref{sse13}, 
we mention some results on quasi-monomial actions 
due to Hoshi, Kang and Kitayama. 
A necessary and sufficient condition for the rationality 
under one-dimensional quasi-monomial actions and 
two-dimensional purely quasi-monomial actions will be given 
via norm residue $2$-symbol. 

In Section \ref{seNB}, 
we treat Noether's problem on rationality. 
Using the unramified Brauer groups, Saltman and Bogomolov were able
to establish counter-examples to Noether's problem 
over algebraically closed field 
for non-abelian $p$-groups of order $p^9$ and $p^6$ respectively.
We mention a result due to Hoshi, Kang and Kunyavskii 
which gives a necessary and sufficient condition for the non-vanishing of the 
unramified Brauer groups for groups of order $p^5$ 
where $p$ is an odd prime number. 

In Section \ref{seAlgTori}, 
we consider the rationality problem for algebraic tori. 
In Subsection \ref{sse31}, 
we recall some results due to Voskresenskii, Kunyavskii, 
Endo and Miyata. 
A birational classification of the algebraic $k$-tori of dimension four 
and five due to Hoshi and Yamasaki will also be given. 
In Subsection \ref{sse32}, we give a detailed account of methods 
related to integral representations of finite groups.

\section{Rationality problem for quasi-monomial actions, generalities}\label{se1}
\subsection{Definitions and examples}\label{sse11}
\begin{definition}
Let $K/k$ and $L/k$ be finitely generated extensions of fields.\\ 
{\rm (1)} $K$ is {\it rational} over
$k$ (for short, $k$-{\it rational}) if $K$ is purely transcendental 
over $k$, i.e. $K\simeq k(x_1,\ldots,x_n)$ 
for some algebraically independent elements $x_1,\ldots,x_n$ over $k$.\\
{\rm (2)} $K$ is {\it stably $k$-rational} if $K(y_1,\ldots,y_m)$ is $k$-rational 
for some $y_1,\ldots,y_m$ such that $y_1,\ldots,y_m$ are
algebraically independent over $K$.\\ 
{\rm (3)} 
$K$ and $L$ are {\it stably $k$-isomorphic} if 
$K(y_1,\ldots,y_m)\simeq L(z_1,\ldots,z_n)$ 
for some algebraically independent elements 
$y_1,\ldots,y_m$ over $K$ and $z_1,\ldots,z_n$ over $L$.\\
{\rm (4)} 
$K$ is 
{\it retract $k$-rational} if there exists a $k$-algebra
$A$ contained in $K$ such that (i) $K$ is the quotient field of
$A$, (ii) there exist a non-zero polynomial $f\in
k[x_1,\ldots,x_n]$ 
and $k$-algebra homomorphisms $\varphi\colon A\to
k[x_1,\ldots,x_n][1/f]$ and $\psi\colon k[x_1,\ldots,x_n][1/f]\to
A$ satisfying $\psi\circ\varphi =1_A$.\\ 
{\rm (5)} $K$ is 
{\it $k$-unirational} 
if $k\subset K\subset k(x_1,\ldots,x_n)$ for some integer $n$. 
\end{definition}
It is not difficult to verify the following implications for an infinite field $k$:
\begin{center}
$k$-rational\ \ $\Rightarrow$\ \ 
stably $k$-rational\ \ $\Rightarrow$\ \ 
retract $k$-rational\ \ $\Rightarrow$\ \ 
$k$-unirational. 
\end{center}

\begin{remark}
In Saltman's original definition of retract $k$-rationality 
(\cite[page 130]{Sal82b}, \cite[Definition 3.1]{Sal84b}), 
a base field $k$ is required to be infinite in order to 
guarantee the existence of sufficiently many $k$-specializations 
(see also Theorem \ref{thgen}).
\end{remark}


\begin{definition}\label{d1.1}
Let $K/k$ be a finite extension of fields and 
$K(x_1,\ldots,x_n)$ be the rational function field over $K$ 
with $n$ variables $x_1,\ldots,x_n$. 
Let $G$ be a finite subgroup of $\fn{Aut}_k(K(x_1,\ldots,x_n))$.\\ 
(1) The action of $G$ on $K(x_1,\ldots,x_n)$ is called {\it quasi-monomial} 
if it satisfies the following three conditions:\\
(i) $\sigma(K)\subset K$ for any $\sigma\in G$;\\
(ii) $K^G=k$, where $K^G$ is the fixed field under the action of $G$;\\
(iii) for any $\sigma\in G$ and any $1 \le j \le n$,
\begin{align}
\sigma(x_j)=c_j(\sigma)\prod_{i=1}^n x_i^{a_{ij}}\label{acts1}
\end{align} 
where
$c_j(\sigma)\in K^\times$ and $[a_{ij}]_{1\le i,j \le n} \in
GL_n(\bm{Z})$.\\
(2) The quasi-monomial action is called {\it purely quasi-monomial} action 
if $c_j(\sigma)=1$ for any $\sigma \in G$ and any $1\le j\le n$ in (iii). \\
(3) The quasi-monomial action is called {\it monomial} action 
if $G$ acts trivially on $K$, i.e.\ $k=K$.\\ 
(4) The quasi-monomial action is called {\it purely monomial} action 
if it is purely quasi-monomial and monomial. 
\end{definition}
We have the following implications: 
\begin{center}
quasi-monomial action\ \ $\Leftarrow$\ \ 
purely quasi-monomial action\\
$\Uparrow$\hspace*{3.5cm}$\Uparrow$~~~~~~~~\\
\ monomial action\ \ $\Leftarrow$\ \ purely monomial action. 
\end{center}

Although there are many variants and results on the rationality problem 
in algebraic geometry and invariant theory, 
we restrict ourselves to the following problem. 

\begin{question}
Let $K/k$ be a finite extension of fields and $G$ be a finite group 
acting on $K(x_1,\ldots,x_n)$ by quasi-monomial $k$-automorphisms. 
Under what situation is the fixed field $K(x_1,\ldots,x_n)^G$ $k$-rational?
\end{question}

This problem includes Noether's problem and the rationality problem for 
algebraic tori (see Example \ref{ex123} below). 
The reader is referred to Swan \cite{Swa83}, 
Manin and Tsfasman \cite{MT86} and 
Colliot-Th\'{e}l\`{e}ne and Sansuc \cite{CTS07} 
for more general survey on the rationality problem, 
and also to Serre \cite{Ser79, Ser02}, 
Knus, Merkurjev, Rost and Tignol \cite{KMRT98}, 
Gille and Szamuely \cite{GS06}
and Berhuy \cite{Ber10} 
for basic tools (e.g. Galois cohomology, Galois descent, Brauer groups, etc.) 
in this area. 

\begin{example}[Typical examples of quasi-monomial actions]\label{ex123}

(1) (Noether's problem). 
When $G$ acts 
on $K(x_1,\ldots,x_n)$ by permutation of the variables 
$x_1,\ldots,x_n$ and 
trivially on $K$, i.e.\ $k=K$, 
the rationality problem of $K(x_1,\ldots,x_n)^G$ 
over $k$ is called Noether's problem. 
We will discuss Noether's problem in 
Section \ref{seNB}.

(2) (Rationality problem for algebraic tori). 
When $G$ acts on $K(x_1,\ldots,x_n)$ 
by purely quasi-monomial $k$-automorphisms and
$G$ is isomorphic to $\fn{Gal}(K/k)$, 
the fixed field $K(x_1,\ldots,x_n)^G$ is a function field of some
algebraic torus defined over $k$ and split over $K$ 
(see Voskresenskii \cite[Chapter 2]{Vos98}).
We will treat 
this in Section \ref{seAlgTori}.

(3) (Rationality problem for Severi-Brauer varieties). 
Assume that $G$ is isomorphic to $\fn{Gal}(K/k)$.  
Take $a_\sigma \in GL_{n+1}(K)$ for each $\sigma \in G$. 
Denote by $\bar{a}_\sigma$ the image of $a_\sigma$ in the
canonical map $GL_{n+1}(K)\to PGL_{n+1}(K)$. Consider the rational
function fields $K(y_0,y_1,\ldots,y_n)$ and $K(x_1$, $\ldots,x_n)$
where $x_i=y_i/y_0$ for $1\le i\le n$. For each $\sigma\in G$,
$a_\sigma$ induces 
a $\fn{Gal}(K/k)$-equivariant automorphism on $K(y_0,y_1,\ldots,y_n)$ 
and $K(x_1,\ldots,x_n)$ 
(note that elements of $K$ in
$K(y_0,\ldots,y_n)$ are acted through $\fn{Gal}(K/k)$). 
Assume furthermore that the map $\gamma:G\to PGL_n(K)$ defined by
$\gamma(\sigma)=\bar{a}_\sigma$ is a 1-cocycle, i.e.\
$\gamma(\sigma \tau)=\gamma(\sigma) \cdot \sigma(\gamma(\tau))$.
Then $G$ induces an action on $K(x_1,\ldots,x_n)$. The fixed field
$K(x_1,\ldots,x_n)^G$ is called a Brauer-field $F_{n,k}(\gamma)$,
i.e.\ the function field of an $n$-dimensional Severi-Brauer
variety over $k$ associated to $\gamma$ 
(see Roquette \cite{Roq63, Roq64}, Kang \cite{Kan90}). 
If $\gamma'$ is a $1$-cocycle which is cohomologous to $\gamma$, it
is easy to see that $F_{n,k}(\gamma) \simeq F_{n,k}(\gamma')$ over
$k$; thus the Brauer-field $F_{n,k}(\gamma)$ depends only on the
cohomology class $[\gamma] \in H^1(G,PGL_n(K))$. 
We can show that a Brauer-field over $k$ is $k$-rational if and only if it is
$k$-unirational (see Serre \cite[page 160]{Ser79}). 
If we assume that each 
column of $a_\sigma$ has precisely one non-zero entry, 
then the action of $G$ on
$K(x_1,\ldots,x_n)$ becomes a quasi-monomial action.
\end{example}

\begin{notation}
Throughout this paper, 
$S_n$ (resp. $A_n$, $D_n$, $C_n$) denotes the symmetric 
(resp. the alternating, the dihedral, the cyclic) group 
of degree $n$ of order $n!$ (resp. $n!/2$, $2n$, $n$). 
\end{notation}

\subsection{Rationality problem for monomial actions}\label{sse12}

Monomial actions 
on $k(x_1,x_2)$ and $k(x_1,x_2,x_3)$ were investigated by Hajja, Kang, Saltman, 
Hoshi, Rikuna, Prokhorov, Kitayama, Yamasaki, etc. 
We list known results for (purely) monomial actions 
on $k(x_1,x_2)$ and $k(x_1,x_2,x_3)$.

\begin{theorem}[{Hajja \cite{Haj87}}]
Let $k$ be a field and $G$ be a finite group acting on $k(x_1,x_2)$
by monomial $k$-automorphisms. Then $k(x_1,x_2)^{G}$ is $k$-rational.
\end{theorem}

\begin{theorem}[{Hajja and Kang \cite{HK92, HK94}, Hoshi and Rikuna \cite{HR08}}]
\label{t1.14}
Let $k$ be a field and $G$ be a finite group acting on
$k(x_1,x_2,x_3)$ by purely monomial $k$-automorphisms. 
Then $k(x_1,x_2,x_3)^G$ is $k$-rational.
\end{theorem}

We define another terminology related to a quasi-monomial action.

\begin{definition}\label{d1.5}
Let $G$ be a finite group acting on $K(x_1,\ldots,x_n)$ by quasi-monomial
$k$-automorphisms with the conditions (i), (ii), (iii) in
Definition \ref{d1.1} (1). Define a group homomorphism $\rho_{\ul{x}}
: G\to GL_n(\bm{Z})$ by $\rho_{\ul{x}}(\sigma)=[a_{ij}]_{1\le
i,j\le n} \in GL_n(\bm{Z})$ for any $\sigma \in G$ where the
matrix $[a_{ij}]_{1\le i,j \le n}$ is given by
$\sigma(x_j)=c_j(\sigma)\prod_{1\le i\le n} x_i^{a_{ij}}$ in (iii)
of Definition \ref{d1.1} (1).
\end{definition}

\begin{proposition}[{\cite[Proposition 1.12]{HKKi}}] \label{p1.6}
Let $G$ be a finite group acting on $K(x_1,\ldots,x_n)$ by
quasi-monomial $k$-automorphisms. Then there exists a normal subgroup
$N$ of $G$ satisfying the following conditions {\rm :} {\rm (i)}
$K(x_1,\ldots,x_n)^N=K^N(y_1,\ldots,y_n)$ where each $y_i$ is of
the form $ax_1^{e_1}x_2^{e_2}\cdots x_n^{e_n}$ with $a\in
K^{\times}$ and $e_i\in \bm{Z}$ {\rm (}we may take $a=1$ if the
action is a purely quasi-monomial action {\rm )}, {\rm (ii)} $G/N$
acts on $K^N(y_1,\ldots,y_n)$ by quasi-monomial $k$-automorphisms,
and {\rm (iii)} $\rho_{\ul{y}}: G/N\to GL_n(\bm{Z})$ is an
injective group homomorphism where
$\rho_{\ul{y}}$ is given as in Definition \ref{d1.5}.
\end{proposition}

Let $G$ be a finite group acting on $k(x_1,x_2,x_3)$ 
by monomial $k$-automorphisms. 
Then we may assume that $G$ is a subgroup of $GL_3(\bm{Z})$ by 
Proposition \ref{p1.6}. 
The rationality problem of $k(x_1,x_2,x_3)^G$ is determined 
by $G$ up to conjugacy in $GL_3(\bm{Z})$. 
There exist $73$ conjugacy classes $(\bm{Z}$-classes) $[G]$ of finite subgroups $G$ 
in $GL_3(\bm{Z})$. 
According to \cite[Table 1]{BBNWZ78}, 
we denote by $[G_{i,j,k}]$ the $k$-th $\bZ$-class 
of the $j$-th $\bQ$-class of the $i$-th crystal system $(1\leq i\leq 7)$ 
of $GL_3(\bZ)$. 
We define 
\begin{align*}
\mathcal{N}:=\bigl\{[G_{1,2,1}],[G_{2,3,1}],[G_{3,1,1}],
[G_{3,3,1}],[G_{4,2,1}],[G_{4,2,2}],[G_{4,3,1}],[G_{4,4,1}]\bigr\}.
\end{align*}

For two groups $G=G_{1,2,1}$, $G_{4,2,2}\in\mathcal{N}$, 
the rationality problem of $k(x_1,x_2,x_3)^G$ over $k$ was 
completely solved by Saltman and Kang: 

\begin{theorem}[Saltman {\cite[Theorem 0.1]{Sal00}, Kang \cite[Theorem 4.4]{Kan05}}]
\label{thSalt} 
\\Let $k$ be a field with {\rm char} $k\neq 2$ and 
$G_{1,2,1}=\langle\sigma\rangle\simeq C_2$ act on 
$k(x_1,x_2,x_3)$ by 
\begin{align*}
\sigma : x_1\ \mapsto\ \tfrac{a_1}{x_1},\ x_2\ \mapsto\ \tfrac{a_2}{x_2},\ 
x_3\ \mapsto\ \tfrac{a_3}{x_3}, 
\quad a_i\in k^\times,\ 1\leq i\leq 3. 
\end{align*}
Then $k(x_1,x_2,x_3)^{G_{1,2,1}}$ is not $k$-rational 
if and only if $[k(\sqrt{a_1},\sqrt{a_2},\sqrt{a_3}):k]=8$. 
If $k(x_1,x_2,x_3)^{G_{1,2,1}}$ is not $k$-rational, then 
it is not retract $k$-rational. 
\end{theorem}
%
\begin{theorem}[Kang {\cite[Theorem 1.8]{Kan04}}]\label{thKan04b}
Let $k$ be a field and $G_{4,2,2}=\langle\sigma\rangle$ $\simeq C_4$ 
act on $k(x_1,x_2,x_3)$ by 
\begin{align*}
\sigma : x_1\ \mapsto\ x_2\ \mapsto\ x_3\ \mapsto\ \tfrac{c}{x_1x_2x_3}\ \mapsto\ x_1,\quad c\in k^\times. 
\end{align*}
Then $k(x_1,x_2,x_3)^{G_{4,2,2}}$ is $k$-rational if and only if 
at least one of the following 
conditions is satisfied: {\rm (i)} {\rm char} $k=2;$ {\rm (ii)} $c\in k^2;$ 
{\rm (iii)} $-4c\in k^4;$ {\rm (iv)} $-1\in k^2$. 
If $k(x_1,x_2,x_3)^{G_{4,2,2}}$ is not $k$-rational, 
then it is not retract $k$-rational. 
\end{theorem}


\begin{theorem}[Hoshi, Kitayama and Yamasaki \cite{HKY11}, Yamasaki \cite{Yam12}]\label{thmainHKY}\\
Let $k$ be a field with {\rm char} $k\neq 2$ 
and $G$ be a finite subgroup of $GL_3(\bZ)$ 
acting on $k(x_1,x_2,x_3)$ by monomial $k$-automorphisms. 

{\rm (1)} {\rm (\cite{HKY11})} 
If $G\not\in \mathcal{N}$, then 
$k(x_1,x_2,x_3)^G$ is $k$-rational except for $G\in [G_{7,1,1}]$;

{\rm (2)} {\rm (\cite{Yam12})} 
If $G\in\mathcal{N}$, then $k(x_1,x_2,x_3)^G$ is not $k$-rational 
for some field $k$ and coefficients $c_j(\sigma)$ 
as in Definition \ref{d1.5}. 
If $k(x_1,x_2,x_3)^G$ is not $k$-rational, 
then it is not retract $k$-rational. 
Indeed, we can obtain the necessary and sufficient condition 
for the $k$-rationality of $k(x_1,x_2,x_3)^G$ in terms of 
$k$ and $c_j(\sigma)$ for each $G\in\mathcal{N}$.
\end{theorem}

For the exceptional case 
$G_{7,1,1}=\langle\tau,\lambda,\sigma\rangle\simeq A_4$, 
the problem can be reduced to the following actions: 
\begin{align*}
\tau &: x_1\ \mapsto\ \tfrac{a}{x_1},\ x_2\ \mapsto\ 
\ep \tfrac{a}{x_2},\ x_3\ \mapsto\ \ep x_3,&
\lambda &: x_1\ \mapsto\ \ep \tfrac{a}{x_1},\ 
x_2\ \mapsto\ \ep x_2,\ x_3\ \mapsto\ \tfrac{a}{x_3},\\
\sigma &: x_1\ \mapsto\  x_2,\ x_2\ \mapsto\  x_3,\ x_3\ \mapsto\ x_1
\end{align*}
where $a\in k^\times$, $\ep=\pm 1$. 
The following gives a partial answer to the problem 
though we do not know 
whether $k(x_1,x_2,x_3)^{G_{7,1,1}}$ is $k$-rational when 
$\ep=-1$ and $[k(\sqrt{a},\sqrt{-1}):k]=4$. 
\begin{theorem}[{\cite[Theorem 1.7]{HKY11}}]\label{thc}
Let $k$ be a field with {\rm char} $k\neq 2$.\\
{\rm (1)} If $\ep=1$, then $k(x_1,x_2,x_3)^{G_{7,1,1}}$ is $k$-rational;\\
{\rm (2)} If $\ep=-1$ and $[k(\sqrt{a},\sqrt{-1}):k]\leq 2$, then 
$k(x_1,x_2,x_3)^{G_{7,1,1}}$ is $k$-rational. 
\end{theorem}

As a consequence of Theorem \ref{thmainHKY} 
and Theorem \ref{thc}, we get the following: 
\begin{theorem}[{\cite[Theorem 1.8]{HKY11}}]\label{thcor}
Let $k$ be a field with {\rm char} $k\neq 2$ and $G$ be a finite group acting on 
$k(x_1,x_2,x_3)$ by monomial $k$-automorphisms. 
Then there exists $L=k(\sqrt{a})$ with $a\in k^\times$ such that 
$L(x_1,x_2,x_3)^G$ is $L$-rational. 
In particular, if $k$ is a quadratically closed field, then $k(x_1,x_2,x_3)^G$ 
is always $k$-rational. 
\end{theorem}

\begin{remark}
Prokhorov \cite[Theorem 5.1]{Pro10} proved 
Theorem \ref{thcor} when $k=L=\bC$ using 
a technique of algebraic geometry (e.g. Segre embedding). 
\end{remark}

\subsection{{Results in \cite{HKKi} for quasi-monomial actions}}\label{sse13}

In this subsection, 
we present some results about quasi-monomial actions 
of dimension $\leq 5$ in Hoshi, Kang and Kitayama \cite{HKKi}. 
It turns out that $K(x_1,\ldots,x_n)^G$ is not even $k$-unirational 
in general. 

\begin{proposition}[{\cite[Proposition 1.13]{HKKi}}] \label{p1.7}
{\rm (1)} Let $G$ be a finite group acting 
on $K(x)$ by purely quasi-monomial $k$-automorphisms.
Then $K(x)^G$ is $k$-rational.

{\rm (2)} Let $G$ be a finite group acting on $K(x)$ by
quasi-monomial $k$-automorphisms. Then $K(x)^G$ is $k$-rational
except for the following case: There exists a normal subgroup
$N$ of $G$ such that {\rm (i)} $G/N=\langle \sigma \rangle \simeq
C_2$, {\rm (ii)} $K(x)^N=k(\alpha)(y)$ with $\alpha^2=a\in
K^{\times}$, $\sigma(\alpha)=-\alpha$ {\rm (}if $\fn{char}k \ne
2${\rm )}, and $\alpha^2+\alpha=a\in K$, $\sigma(\alpha)=\alpha+1$
{\rm (}if $\fn{char} k=2${\rm )}, {\rm (iii)} $\sigma\cdot y=b/y$
for some $b\in k^{\times}$.

For the exceptional case, $K(x)^G=k(\alpha) (y)^{G/N}$ is
$k$-rational if and only if the norm residue $2$-symbol
$(a,b)_k=0$ {\rm (}if $\fn{char}k\ne 2${\rm )}, and $[a,b)_k=0$ 
{\rm (}if $\fn{char}k=2${\rm )}.

Moreover, if $K(x)^G$ is not $k$-rational, then the Brauer group 
$\fn{Br}(k)$ is not trivial and $K(x)^G$ is not $k$-unirational.
\end{proposition}

For the definition of the norm residue $2$-symbols $(a,b)_k$ and
$[a,b)_k$, see \cite[Chapter 11]{Dra83}.

\begin{theorem}[{\cite[Theorem 1.14]{HKKi}}] \label{t1.8}
Let $G$ be a finite group and $G$ act on $K(x,y)$ by purely
quasi-monomial $k$-automorphisms. Define 
\[
N=\{\sigma\in G\mid\sigma(x)=x,~ \sigma(y)=y\},\ 
H=\{\sigma\in G\mid\sigma(\alpha)=\alpha\ {\rm for\ all}\ \alpha\in K\}.
\]
Then $K(x,y)^G$ is $k$-rational except possibly for the following situation: 
{\rm (1)} $\fn{char}k\ne 2$ and {\rm (2)}
$(G/N,HN/N)\simeq (C_4,C_2)$ or $(D_4,C_2)$.

More precisely, in the exceptional situation we may choose $u,v\in
k(x,y)$ satisfying that $k(x,y)^{HN/N}=k(u,v)$ 
{\rm (}and therefore $K(x,y)^{HN/N}=K(u,v)${\rm )} such that
\begin{enumerate}
\item[{\rm (i)}] when $(G/N,HN/N)\simeq (C_4,C_2)$,
$K^N=k(\sqrt{a})$ for some $a\in k\backslash k^2$, $G/N=\langle
\sigma \rangle \simeq C_4$, then $\sigma$ acts on $K^N(u,v)$ by
$\sigma: \sqrt{a} \mapsto -\sqrt{a}$, $u\mapsto \frac{1}{u}$,
$v\mapsto -\frac{1}{v}$; or \item[{\rm (ii)}] when
$(G/N,HN/N)\simeq (D_4,C_2)$, $K^N=k(\sqrt{a},\sqrt{b})$ is a
biquadratic extension of $k$ with $a,b\in k\backslash k^2$,
$G/N=\langle \sigma,\tau \rangle \simeq D_4$, then $\sigma$ and
$\tau$ act on $K^N(u,v)$ by $\sigma:\sqrt{a} \mapsto -\sqrt{a}$,
$\sqrt{b}\mapsto \sqrt{b}$, $u\mapsto\frac{1}{u}$, $v\mapsto
-\frac{1}{v}$, $\tau: \sqrt{a}\mapsto \sqrt{a}$, $\sqrt{b}\mapsto
-\sqrt{b}$, $u\mapsto u$, $v\mapsto -v$.
\end{enumerate}
For case {\rm (i)}, $K(x,y)^G$ is $k$-rational if and only if the
norm residue $2$-symbol $(a,-1)_k=0$. For case {\rm (ii)},
$K(x,y)^G$ is $k$-rational if and only if $(a,-b)_k=0$.

Moreover, if $K(x,y)^G$ is not $k$-rational, then the Brauer group 
$\fn{Br}(k)$ is not trivial and $K(x,y)^G$ is not $k$-unirational.
\end{theorem}


Saltman \cite[Section 3]{Sal90a} discussed the relationship of $K(M)^G$ and the
embedding problem in Galois theory. 
We reformulate the two exceptional cases of 
Theorem \ref{t1.8} in terms of the embedding problem.

\begin{proposition}[{\cite[Proposition 4.3]{HKKi}}] \label{p4.2}
{\rm (1)}
Let $k$ be a field with $\fn{char}k\ne 2$, $a\in k\backslash k^2$ and $K=k(\sqrt{a})$.
Let $G=\langle \sigma\rangle$ act on the rational function field $K(u,v)$ by
\[
\sigma: \sqrt{a} \mapsto -\sqrt{a},~u\mapsto \tfrac{1}{u},~ v\mapsto -\tfrac{1}{v}.
\]
Then $K(u,v)^G$ is $k$-rational if and only if 
there exists a Galois extension $L/k$ such that 
$k\subset k(\sqrt{a})\subset L$ and $\fn{Gal}(L/k)$ is isomorphic to $C_4$.

{\rm (2)}
Let $k$ be a field with $\fn{char}k\ne 2$ and 
$a,b\in k\backslash k^2$ such that $[K:k]=4$ where $K=k(\sqrt{a},\sqrt{b})$.
Let $G=\langle \sigma,\tau \rangle$ act on the rational function field $K(u,v)$ by
\begin{align*}
\sigma &: \sqrt{a}\mapsto -\sqrt{a},~ \sqrt{b}\mapsto\sqrt{b},~ u\mapsto \tfrac{1}{u},~v\mapsto -\tfrac{1}{v}, \\
\tau &: \sqrt{a}\mapsto \sqrt{a},~ \sqrt{b}\mapsto -\sqrt{b},~ u\mapsto u,~ v\mapsto -v.
\end{align*}
Then $K(u,v)^G$ is $k$-rational if and only if 
there exists a Galois extension $L/k$ such that 
$k\subset K \subset L$ and $\fn{Gal}(L/k)$ is isomorphic to $D_4$.
\end{proposition}


{}From the action of $G/N=\langle\sigma,\tau\rangle\simeq D_4$ on $K(x,y)$ 
in the exceptional case (ii) of Theorem \ref{t1.8}, 
we obtain the following example. 

\begin{example}[{\cite[Example 4.4]{HKKi}}] \label{ex4.3}
Let $k$ be a field with $\fn{char}k\ne 2$ and 
$K=k(\sqrt{a},\sqrt{b})$ be a biquadratic extension of $k$. 
We consider the following actions of $k$-automorphisms of 
$D_4=\langle\sigma,\tau \rangle$ on $K(x,y)$:
\begin{align*}
\sigma_a &: \sqrt{a} \mapsto -\sqrt{a},~ \sqrt{b}\mapsto \sqrt{b},~ x\mapsto y,~ y\mapsto\tfrac{1}{x}, \\
\sigma_b &: \sqrt{a}\mapsto \sqrt{a},~ \sqrt{b}\mapsto -\sqrt{b},~ x\mapsto y,~ y\mapsto \tfrac{1}{x}, \\
\sigma_{ab} &: \sqrt{a} \mapsto -\sqrt{a},~ \sqrt{b}\mapsto -\sqrt{b},~ x\mapsto y,~ y\mapsto \tfrac{1}{x}, \\
\tau_a &: \sqrt{a} \mapsto -\sqrt{a},~ \sqrt{b}\mapsto \sqrt{b},~ x\mapsto y,~ y\mapsto x, \\
\tau_b &: \sqrt{a}\mapsto \sqrt{a},~ \sqrt{b}\mapsto -\sqrt{b},~ x\mapsto y,~ y\mapsto x, \\
\tau_{ab} &: \sqrt{a}\mapsto -\sqrt{a},~ \sqrt{b}\mapsto -\sqrt{b},~ x\mapsto y,~ y\mapsto x.
\end{align*}
Define $L_{a,b}=K(x,y)^{\langle\sigma_a,\tau_b\rangle}$. Then\\
(1) $L_{a,b}$ is $k$-rational $\Leftrightarrow$ 
$L_{a,ab}$ is $k$-rational $\Leftrightarrow$ $(a,-b)_k=0$;\\
(2) $L_{b,a}$ is $k$-rational $\Leftrightarrow$ 
$L_{b,ab}$ is $k$-rational $\Leftrightarrow$ $(b,-a)_k=0$;\\
(3) $L_{ab,a}$ is $k$-rational $\Leftrightarrow$ 
$L_{ab,b}$ is $k$-rational $\Leftrightarrow$ $(a,b)_k=0$.

In particular, if $\sqrt{-1}\in k$, then the obstructions to the
rationality of the above fixed fields over $k$ coincide, i.e. they
are reduced to the same condition, $(a,b)_k=0$.

On the other hand, consider the case $k=\bm{Q}$,
$K=\bm{Q}(\sqrt{-1},\sqrt{p})$ where
$p$ is a prime number with $p\equiv 1$ (mod 4). Then
$L_{a,b}$ and $L_{a,ab}$ are not $\bm{Q}$-rational 
because $(-1,-p)_{\bm{Q}}=(-1,-1)_{\bm{Q}} \neq 0$,
while $L_{b,a}$, $L_{b,ab}$,
$L_{ab,a}$ and $L_{ab,b}$ are $\bm{Q}$-rational.
\end{example}


The following gives an equivalent definition of 
purely quasi-monomial actions. 

\begin{definition} \label{d1.9}
Let $G$ be a finite group. 
A {\it $G$-lattice} $M$ is a finitely 
generated $\bm{Z}[G]$-module which is $\bm{Z}$-free as an abelian
group, i.e.\ $M=\bigoplus_{1\le i\le n} \bm{Z}\cdot x_i$ with a
$\bm{Z}[G]$-module structure. 
Let $K/k$ be a field extension such
that $G$ acts on $K$ with $K^G=k$. 
We define a purely quasi-monomial action of $G$ on the rational function field
$K(x_1,\ldots,x_n)$ by 
$\sigma\cdot x_j=c_j(\sigma)\prod_{1\le i\le n} x_i^{a_{ij}} \in 
K(x_1,\ldots,x_n)$ 
when $\sigma\cdot x_j=\sum_{1\le i\le n} a_{ij} x_i \in M$ 
and its fixed field is denoted by $K(M)^G$.
\end{definition}

With the aid of Theorem \ref{t1.8}, we are able to show
that $k(M)^G$ is $k$-rational whenever $M$ is a decomposable
$G$-lattice of $\bm{Z}$-rank 4. 
There are $710$ $G$-lattices of $\bm{Z}$-rank 4 and the total number
of decomposable ones is 415. 

\begin{theorem}[{\cite[Theorem 1.16]{HKKi}}] \label{t1.10}
Let $k$ be a field, $G$ be a finite group and $M$ be a $G$-lattice
with $\fn{rank}_{\bm{Z}} M=4$ such that $G$ acts on $k(M)$ by
purely monomial $k$-automorphisms. If $M$ is decomposable, i.e.\
$M=M_1\oplus M_2$ as $\bm{Z}[G]$-modules where $1\le
\fn{rank}_{\bm{Z}} M_1 \le 3$, then $k(M)^G$ is $k$-rational.
\end{theorem}


\begin{proposition}[{\cite[Proposition 5.2]{HKKi}}] \label{p5.2}
Let $G$ be a finite group and $M$ be a $G$-lattice. 
Assume that $M=M_1\oplus M_2\oplus M_3$ as $\bm{Z}[G]$-modules where
$\fn{rank}_{\bm{Z}} M_1=\fn{rank}_{\bm{Z}} M_2=2$ and
$\fn{rank}_{\bm{Z}}M_3=1$. 
Let $k$ be a field, and let $G$ act on $k(M)$ 
by purely monomial $k$-automorphisms. 
Then $k(M)^G$ is $k$-rational.
\end{proposition}

\begin{theorem}[{\cite[Theorem 6.2, Theorem 6.4]{HKKi}}] \label{t6.2}
Let $G$ be a finite group and $M$ be a $G$-lattice. 
Assume that {\rm (i)} $M=M_1\oplus M_2$ as $\bm{Z}[G]$-modules where
$\fn{rank}_{\bm{Z}}M_1=3$ and $\fn{rank}_{\bm{Z}}M_2=2$, 
{\rm (ii)} either $M_1$ or $M_2$ is a faithful $G$-lattice. 
Let $k$ be a field, and let $G$ act on $k(M)$ by purely monomial
$k$-automorphisms. Then $k(M)^G$ is $k$-rational except the
following situation: $\fn{char}k\ne 2$, $G=\langle
\sigma,\tau\rangle \simeq D_4$ and $M_1=\bigoplus_{1\le i\le 3}
\bm{Z} x_i$, $M_2=\bigoplus_{1\le j\le 2} \bm{Z}y_j$ such that
$\sigma:x_1\leftrightarrow x_2$, $x_3\mapsto -x_1-x_2-x_3$,
$y_1\mapsto y_2\mapsto -y_1$, $\tau: x_1\leftrightarrow x_3$,
$x_2\mapsto -x_1-x_2-x_3$, $y_1\leftrightarrow y_2$ where the
$\bm{Z}[G]$-module structure of $M$ is written additively. 

For the exceptional case, $k(M)^G$ is not retract $k$-rational. 
In particular, if $G=\langle\sigma,\tau\rangle\simeq D_4$ acts on the rational 
function field $k(x_1,x_2,x_3,x_4,x_5)$ by $k$-automorphisms 
\begin{align*}
\sigma &: x_1\mapsto x_2,~ x_2\mapsto x_1,~ x_3\mapsto \tfrac{1}{x_1x_2x_3},~ x_4\mapsto x_5,~ x_5\mapsto \tfrac{1}{x_4}, \\
\tau &: x_1\mapsto x_3,~ x_2\mapsto \tfrac{1}{x_1x_2x_3},~ x_3\mapsto x_1,~ x_4\mapsto x_5,~ x_5\mapsto x_4,
\end{align*}
then $k(x_1,x_2,x_3,x_4,x_5)^G$ is not retract $k$-rational.
\end{theorem}

\begin{remark}
The exceptional case of Theorem \ref{t6.2} gives an example of 
purely monomial action of $G\simeq D_4$ whose invariant field $k(M)^G$ 
is not retract $k$-rational even over an algebraically closed field $k$ 
(cf. Theorem \ref{t1.14} and Theorem \ref{thcor}).
\end{remark}

We can deduce the exceptional case of Theorem \ref{t6.2} from the 
following theorem, where the action is not even quasi-monomial.
\begin{theorem}[{\cite[Theorem 6.3]{HKKi}}] \label{t6.3}
Let $k$ be a field with $\fn{char}k\ne 2$ and
$G=\langle\rho\rangle \simeq C_2$ act on 
$k(x_1,x_2,x_3,x_4)$ by $k$-automorphisms defined as
\[
\rho: x_1\mapsto -x_1,~ x_2\mapsto \tfrac{x_4}{x_2},~ x_3\mapsto 
\tfrac{(x_4-1)(x_4-x_1^2)}{x_3},~ x_4\mapsto x_4.
\]
Then $k(x_1,x_2,x_3,x_4)^G$ is not retract $k$-rational.
\end{theorem}

\section{Noether's problem and unramified Brauer group}\label{seNB}

Let $G$ be a finite group acting on the
rational function field $k(x_g\mid g\in G)$ by $k$-automorphisms so that
$g\cdot x_h=x_{gh}$ for any $g,h\in G$. Denote by $k(G)$ the fixed
field $k(x_g\mid g\in G)^G$. Noether's problem asks whether $k(G)$ is
rational (= purely transcendental) over $k$. It is related to the
inverse Galois problem, to the existence of generic $G$-Galois
extensions over $k$, and to the existence of versal $G$-torsors over
$k$-rational field extensions 
(see Garibaldi, Merkurjev and Serre \cite[33.1, page 86]{GMS03}). 
%
For example, Saltman proved the following theorem (see also \cite[Chapter 5]{JLY02}):
\begin{theorem}[see Saltman \cite{Sal82a, Sal82b} for details]\label{thgen}
Let $G$ be a finite group. Assume that $k$ is an infinite field. 
Then the following conditions are equivalent:\\
{\rm (1)} $k(G)$ is retract $k$-rational;\\
{\rm (2)} $G$ has the lifting property over $k$;\\
{\rm (3)} there exists a generic $G$-Galois extension 
$($resp. generic $G$-polynomial$)$ over $k$.
\end{theorem}

We recall some known results on Noether's problem. 

\begin{theorem}[{Fischer \cite{Fis15}, see also Swan \cite[Theorem 6.1]{Swa83}}]\label{thFis}
Let $G$ be a finite abelian group with exponent $e$. 
Assume that {\rm (i)} either {\rm char} $k=0$ or {\rm char} $k>0$ with 
{\rm char} $k$ $\not{|}$ $e$, and 
{\rm (ii)} $k$ contains a primitive $e$-th root of unity. 
Then $k(G)$ is $k$-rational. 
\end{theorem}

Kuniyoshi established the following theorem for $p$-groups 
(see \cite{Kun56}, Proceedings of the international symposium 
on algebraic number theory, Tokyo \& Nikko, 1955).
\begin{theorem}[Kuniyoshi \cite{Kun54, Kun55, Kun56}]
Let $G$ be a $p$-group and $k$ be a field with {\rm char} $k=p>0$. 
Then $k(G)$ is $k$-rational. 
\end{theorem}

Swan \cite{Swa69} showed that 
$\bQ(C_{47})$ is not $\bQ$-rational using 
Masuda's idea \cite{Mas55, Mas68}. 
This is the first negative example to Noether's problem. 
After efforts of many mathematicians  (e.g. Voskresenskii, Endo and Miyata), 
Noether's problem for abelian groups was solved by Lenstra \cite{Len74}. 
The reader is referred to Swan's paper \cite{Swa83} for a survey of this problem.


On the other hand, just a handful of results about Noether's problem
are obtained when the groups are not abelian. 

\begin{theorem}[Chu and Kang \cite{CK01}] \label{t1.6}
Let $p$ be any prime number and $G$ be a $p$-group of order $\le p^4$
and of exponent $e$. If $k$ is a field containing a primitive $e$-th root
of unity, then $k(G)$ is $k$-rational.
\end{theorem}

\begin{theorem}[{Serre \cite[Chapter IX]{GMS03}}]
Let $G$ be a group with a $2$-Sylow subgroup which is cyclic of 
order $\geq 8$ or the generalized quaternion $Q_{16}$ of order $16$. 
Then $\bQ(G)$ is not $\bQ$-rational. 
\end{theorem}

\begin{theorem}[Chu, Hu, Kang and Prokhorov \cite{CHKP08}] \label{t1.7}
Let $G$ be a group of order $2^5$ and of exponent $e$. 
If $k$ is a field containing a primitive $e$-th root
of unity, then $k(G)$ is $k$-rational.
\end{theorem}

The notion of the unramified Brauer group of $K$ over $k$, denoted by
$\fn{Br}_{v,k}(K)$, was introduced by Saltman \cite{Sal84a}.
\begin{definition}[{Saltman \cite[Definition 3.1]{Sal84a}, \cite[page 56]{Sal85}}]\label{d1.2}
Let $K/k$ be an extension of fields. 
{\it The unramified Brauer group} $\fn{Br}_{v,k}(K)$ of $K$ over $k$ is 
defined by 
$\fn{Br}_{v,k}(K)=\bigcap_R \fn{Image} \{ \fn{Br}(R)\to
\fn{Br}(K)\}$ where $\fn{Br}(R)\to \fn{Br}(K)$ is the natural map of
Brauer groups and $R$ runs over all the discrete valuation rings $R$
such that $k\subset R\subset K$ and $K$ is the quotient field of
$R$.
\end{definition}

\begin{lemma}[{Saltman \cite{Sal84a}, \cite[Proposition 1.8]{Sal85}, \cite{Sal87}}] \label{l1.3}
If 
$K$ is retract $k$-rational, then
the natural map $\fn{Br}(k)\to \fn{Br}_{v,k} (K)$ is an isomorphism.
In particular, if $k$ is an algebraically closed field and $K$ is
retract $k$-rational, then $\fn{Br}_{v,k}(K)=0$.
\end{lemma}

\begin{theorem}[{Bogomolov \cite[Theorem 3.1]{Bog88}, 
Saltman \cite[Theorem 12]{Sal90b}}] \label{t1.4}
Let $G$ be a finite group and $k$ be an algebraically closed field with 
$\gcd \{|G|,\fn{char}k\}=1$. 
Then $\fn{Br}_{v,k}(k(G))$ is isomorphic to the group $B_0(G)$ defined by
\[
B_0(G)=\bigcap_A \fn{Ker}\{\fn{res}: H^2(G,\bm{Q}/\bm{Z})\to H^2(A,\bm{Q}/\bm{Z})\}
\]
where $A$ runs over all the bicyclic subgroups of $G$ $($a group $A$
is called bicyclic if $A$ is either a cyclic group or a direct
product of two cyclic groups$)$.
\end{theorem}

We may call $B_0(G)$ the Bogomolov multiplier of $G$ 
since $H^2(G,\bm{Q}/\bm{Z})$ is isomorphic to the Schur multiplier 
$H_2(G,\bm{Z})$ of $G$ (see Karpilovsky \cite{Kar87}, Kunyavskii \cite{Kun10}). 
Because of Theorem \ref{t1.4}, we will not 
distinguish $B_0(G)$ and $\fn{Br}_{v,k}(k(G))$ when $k$ is
algebraically closed and $\gcd \{|G|,\fn{char}k\}=1$. 

Using the unramified Brauer groups, Saltman and Bogomolov are able
to establish counter-examples to Noether's problem for non-abelian
$p$-groups over algebraically closed field.

\begin{theorem}[Saltman, Bogomolov] \label{t1.5}
Let $p$ be any prime number and $k$ be any algebraically closed field
with $\fn{char}k\ne p$.

{\rm (1) (Saltman \cite{Sal84a})} There exists a group $G$ of order $p^9$
such that $B_0(G)\ne 0$. In particular, $k(G)$ is not retract
$k$-rational. Thus $k(G)$ is not $k$-rational.

{\rm (2) (Bogomolov \cite{Bog88})} There exists a group $G$ of order
$p^6$ such that $B_0(G)\ne 0$. Thus $k(G)$ is not $($retract$)$ $k$-rational.
\end{theorem}


\begin{example} \label{tt1.8}
{\rm (1)} Let $p$ be an odd prime number. 
If $G$ is $p$-group of order $\leq p^4$ or $2$-group of order $\leq 2^5$, 
then $B_0(G)=0$ (see Theorem \ref{t1.6} and Theorem \ref{t1.7}). 

{\rm (2)} Working on $p$-groups, 
Bogomolov \cite{Bog88} developed a lot of techniques 
and interesting results on $B_0(G)$. 
Bogomolov claimed that

{\rm (i) (\cite[Lemma 4.11]{Bog88})} 
If $G$ is a $p$-group with
$B_0(G)\ne 0$ and $G/[G,G]\simeq C_p\times C_p$, then $p\ge 5$ and
$|G|\geq p^7;$

{\rm (ii) (\cite[Lemma 5.6]{Bog88})} 
If $G$ is a $p$-group of order $\le p^5$, then $B_0(G)=0$.

Because of {\rm (ii)}, Bogomolov proposed to
classify all the groups $G$ with $|G|=p^6$ satisfying $B_0(G)\ne 0$
(\cite[Remark 1, page 479]{Bog88}).  
(However, it turns out that (i) and (ii) are not correct, see 
Theorem \ref{t1.9} and Theorem \ref{t1.11}.)
\end{example}

There are 267 non-isomorphic groups of order $2^6$ 
(see Hall and Senior \cite{HS64}). 

\begin{theorem}[Chu, Hu, Kang and Kunyavskii \cite{CHKK10}] \label{thCHKK10}
Let $G=G(2^6,i)$ be the $i$-th group of order $2^6$ 
in the database of GAP {\rm \cite{GAP}} $(1\leq i\leq 267)$.

{\rm (1)} $B_0(G)\ne 0$ if and only if $G=G(2^6,i)$ where $149\le i\le 151$,
$170\le i\le 172$, $177\le i\le 178$, or $i=182$.

{\rm (2)} If $B_0(G)= 0$ and $k$ is an algebraically closed field
with $\fn{char}k\ne 2$, then $k(G)$ is $k$-rational except
possibly for groups $G=G(2^6,i)$ with $241\le i\le 245$.
\end{theorem}

It came as a surprise that Moravec \cite{Mor12}
disproved Example \ref{tt1.8} (2) (i), (ii).

\begin{theorem}[{Moravec \cite[Section 5]{Mor12}}] \label{t1.9}
If $G$ is a group of order $3^5$, then $B_0(G)\ne 0$ if and only if
$G=G(3^5,i)$ with $28\le i\le 30$, where $G(3^5,i)$ is the $i$-th
group among groups of order $3^5$ in the database of GAP {\rm \cite{GAP}}.
\end{theorem}

There exist $67$ non-isomorphic groups of order $3^5$. 
Moravec proves Theorem \ref{t1.9} by using computer calculations 
(GAP functions for computing $B_0(G)$ are available at his website 
\verb+www.fmf.uni-1j.si/~moravec/b0g.g+). 

\begin{definition}\label{d1.10}
Two $p$-groups $G_1$ and $G_2$ are called {\it isoclinic} if there exist
group isomorphisms $\theta\colon G_1/Z(G_1) \to G_2/Z(G_2)$ and
$\phi\colon [G_1,G_1]\to [G_2,G_2]$ such that $\phi([g,h])$
$=[g',h']$ for any $g,h\in G_1$ with $g'\in \theta(gZ(G_1))$, $h'\in
\theta(hZ(G_1))$.

For a prime number $p$ and a fixed integer $n$, let $G_n(p)$ be the
set of all non-isomorphic groups of order $p^n$. In $G_n(p)$
consider an equivalence relation: two groups $G_1$ and $G_2$ are
equivalent if and only if they are isoclinic. Each equivalence class
of $G_n(p)$ is called an {\it isoclinism family}.
\end{definition}


There exist  
\[
2p+61+\gcd \{4,p-1\}+2\gcd\{3,p-1\}
\]
non-isomorphic groups of order $p^5$ $(p\geq 5)$. 
For $p\geq 3$, there are precisely 10 isoclinism families $\Phi_1,\ldots\Phi_{10}$ 
for groups of order $p^5$ (see \cite[pages 619--621]{Jam80}). 

Hoshi, Kang and Kunyavskii \cite{HKKu} 
reached to the following theorem which asserts that 
the non-vanishing of $B_0(G)$ is determined by 
the isoclinism family of $G$ when $G$ is of order $p^5$ . 

\begin{theorem}[{\cite[Theorem 1.12]{HKKu}}] \label{t1.11}
Let $p$ be any odd prime number and $G$ be a group of order $p^5$. 
Then $B_0(G)\ne 0$ if and only if $G$ belongs to the isoclinism family
$\Phi_{10}$. Each group $G$ in the family $\Phi_{10}$ satisfies the
condition $G/[G,G] \simeq C_p\times C_p$. There are precisely $3$
groups in this family if $p=3$. For $p\ge 5$, the total number of
non-isomorphic groups in this family is
\[
1+\gcd\{4,p-1\}+\gcd \{3,p-1\}.
\]
\end{theorem}

\begin{remark}
For groups of order $2^6$, there exist $27$ isoclinism families 
(see \cite[page~147]{JNOB90}). 
We see that Theorem \ref{thCHKK10} can be rephrased 
in terms of the isoclinism family as follows. 
Let $G$ be a group of order $2^6$. 

{\rm (1)} $B_0(G)\ne 0$ if and only if $G$ belongs to the $16$th
isoclinism family $\Phi_{16}$;

{\rm (2)} If $B_0(G)= 0$ and $k$ is an algebraically closed field
with $\fn{char}k\ne 2$, then $k(G)$ is $k$-rational except
possibly for groups $G$ belonging to the $13$th isoclinism family $\Phi_{13}$.
\end{remark}

By Theorem \ref{thCHKK10} and Theorem \ref{t1.11}, 
we get the following theorem (this theorem is optimal by 
Theorem \ref{t1.6}, Theorem \ref{t1.7} and Lemma \ref{l1.3}).

\begin{theorem}[{\cite[Theorem 1.13]{HKKu}}] \label{t1.12}
If $2^6\mid n$ or $p^5\mid n$ for some 
odd prime number $p$, then there exists a group $G$ of order $n$ such
that $B_0(G)\ne 0$. 
\end{theorem}

The following result supplements Moravec's
result (Theorem \ref{t1.9}).

\begin{theorem}[Chu, Hoshi, Hu and Kang \cite{CHHK}] \label{t1.15}
Let $G$ be a group of order $3^5$ and of exponent $e$. If $k$ is a
field containing a primitive $e$-th root of unity and $B_0(G)=0$,
then $k(G)$ is $k$-rational except possibly for groups $G \in
\Phi_7$, i.e.\ $G=G(3^5,i)$ with $56 \le i\le 60$.
\end{theorem}


{}From observations in the proof of Theorem \ref{t1.11}, 
we raised the following question. 

\begin{question}[{\cite[Question 1.11]{HKKu}}] \label{quest-iso}
Let $G_1$ and $G_2$ be isoclinic $p$-groups. Is it true that the
fields $k(G_1)$ and $k(G_2)$ are stably $k$-isomorphic?
\end{question}

Recently, Moravec \cite{Mor} (arXiv:1203.2422) announced that 
if $G_1$ and $G_2$ are isoclinic, 
then $B_0(G_1)\simeq B_0(G_2)$. 
Furthermore, 
Bogomolov and B\"ohning \cite[Theorem 3.2]{BB} (arXiv:1204.4747) 
announced that 
the answer to Question \ref{quest-iso} is affirmative. 


\section{Rationality problem for algebraic tori}\label{seAlgTori}

\subsection{Low-dimensional cases}\label{sse31}

Let $L$ be a finite Galois extension of $k$ and $G={\rm Gal}(L/k)$ 
be the Galois group of the extension $L/k$. 
Let $M=\bigoplus_{1\leq i\leq n}\bZ\cdot u_i$ be a $G$-lattice with 
a $\bZ$-basis $\{u_1,\ldots,u_n\}$, 
i.e. finitely generated $\bZ[G]$-module 
which is $\bZ$-free as an abelian group. 
Let $G$ act on the rational function field $L(x_1,\ldots,x_n)$ 
over $L$ with $n$ variables $x_1,\ldots,x_n$ by 
\begin{align}
\sigma(x_i)=\prod_{j=1}^n x_j^{a_{ij}},\quad 1\leq i\leq n\label{acts}
\end{align}
for any $\sigma\in G$, when $\sigma (u_i)=\sum_{j=1}^n a_{ij} u_j$, 
$a_{ij}\in\bZ$. 
The field $L(x_1,\ldots,x_n)$ with this action of $G$ will be denoted 
by $L(M)$. 

Note that we changed the definition of the action of $\sigma\in G$ 
on $L(x_1,\ldots,x_n)$ from (\ref{acts1}) by the following reason: 

The category of $G$-lattices is anti-equivalent to the 
category of algebraic $k$-tori which split over $L$ 
(see \cite[page 27, Example 6]{Vos98}, \cite[Proposition 20.17]{KMRT98}). 
Indeed, if $T$ is an algebraic $k$-torus, then the character 
group $X(T)={\rm Hom}(T,\bG_m)$ of $T$ becomes $G$-lattice. 
Conversely, for a $G$-lattice $M$, there exists an algebraic $k$-torus 
$T$ which splits over $L$ such that $X(T)$ is isomorphic to $M$ as a $G$-lattice. 

The invariant field $L(M)^G$ 
may be identified with the function field of 
$T$ and is 
always $k$-unirational 
(see \cite[page 40, Example 21]{Vos98}).
Tori of dimension $n$ over $k$ correspond bijectively 
to the elements of the set $H^1(\mathcal{G},GL_n(\bZ))$ 
via Galois descent 
where $\mathcal{G}={\rm Gal}(k_{\rm s}/k)$ since 
${\rm Aut}(\bG_m^n)=GL_n(\bZ)$. 
The $k$-torus $T$ of dimension $n$ is determined uniquely by the integral 
representation $h : \mathcal{G}\rightarrow GL_n(\bZ)$ up to conjugacy, 
and $h(\mathcal{G})$ is a finite subgroup of $GL_n(\bZ)$ 
(see \cite[page 57, Section 4.9]{Vos98}). 

Let $K/k$ be a separable field extension of degree $n$ 
and $L/k$ be the Galois closure of $K/k$. 
Let $G={\rm Gal}(L/k)$ and $H={\rm Gal}(L/K)$. 
The Galois group $G$ may be regarded as a transitive subgroup of $S_n$. 
We have an exact sequence of $\bZ[G]$-modules 
\[
0\longrightarrow I_{G/H}\longrightarrow \bZ[G/H]
\overset{\ep}{\longrightarrow} \bZ\longrightarrow 0,
\]
where $\ep : \bZ[G/H]\rightarrow \bZ$ is the augmentation map. 
Taking the duals, we get an exact sequence of $\bZ[G]$-modules
\[
0\longrightarrow \bZ\longrightarrow \bZ[G/H]
\longrightarrow J_{G/H}\longrightarrow 0.
\]
In particular, $J_{G/H}={\rm Hom}_\bZ(I_{G/H},\bZ)$ 
(called the Chevalley module) is of rank $n-1$. 
Furthermore, we get an exact sequence of algebraic $k$-tori
\[
1\longrightarrow R_{K/k}^{(1)}(\bG_{m,K})\longrightarrow R_{K/k}(\bG_{m,K})
\overset{N_{K/k}}{\longrightarrow} \bG_m\longrightarrow 1.
\]
Here $R_{K/k}(\bG_{m,K})$ is the Weil restriction of 
$\bG_{m,K}$ by the extension $K/k$ and 
$R_{K/k}^{(1)}(\bG_{m,K})$ is the norm one torus of $K/k$ 
whose character group is $J_{G/H}$.

Write $J_{G/H}=\oplus_{1\leq i\leq n-1}\bZ x_i$. 
Then the action of $G$ on $L(J_{G/H})=L(x_1,\ldots,x_{n-1})$ is 
nothing but (\ref{acts}). 

All the $1$-dimensional algebraic $k$-tori $T$, 
i.e. the trivial torus $\bG_m$ and the norm one 
torus $R_{L/k}^{(1)}(\bG_m)$ with $[L:k]=2$, are $k$-rational. 
A birational classification of the 2-dimensional
algebraic tori and the 3-dimensional algebraic tori 
was given 
by Voskresenskii \cite{Vos67} and Kunyavskii \cite{Kun90} respectively. 

\begin{theorem}[Voskresenskii, Kunyavskii]\label{tht23}
Let $k$ be a field.\\
{\rm (1) (Voskresenskii \cite{Vos67})} 
All the two-dimensional algebraic
$k$-tori are $k$-rational. In particular, $K(x_1,x_2)^G$ is
always $k$-rational if $G$ is isomorphic to ${\rm Gal}(K/k)$ and $G$ acts on
$K(x_1,x_2)$ by purely quasi-monomial $k$-automorphisms.\\
{\rm (2) (Kunyavskii \cite{Kun90}, see also Kang \cite[Section 1]{Kan12} for the last statement)} 
All the three-dimensional algebraic
$k$-tori are $k$-rational except for the $15$ cases in the
list of {\rm \cite[Theorem 1]{Kun90}}. For the exceptional $15$
cases, they are not $k$-rational; in fact, they are even not
retract $k$-rational.
\end{theorem}

For $4$-dimensional case, some birational invariants are computed by 
Popov \cite{Pop98}. 


We list known results about the rationality problem for 
norm one tori $R^{(1)}_{K/k}(\bG_m)$ of $K/k$. 

When an extension $K/k$ is Galois, the following theorem is known. 

\begin{theorem}\label{th13}
Let $K/k$ be a finite Galois field extension and $G={\rm Gal}(K/k)$.\\
{\rm (1)} {\rm (Endo and Miyata \cite[Theorem 1.5]{EM74}, 
Saltman \cite[Theorem 3.14]{Sal84b})}\\ 
$R^{(1)}_{K/k}(\bG_m)$ is retract $k$-rational 
if and only if all the Sylow subgroups of $G$ are cyclic.\\
{\rm (2)} {\rm (Endo and Miyata \cite[Theorem 2.3]{EM74})} 
$R^{(1)}_{K/k}(\bG_m)$ is stably $k$-rational 
if and only if $G=C_m$ or 
$G=C_n\times \langle\sigma,\tau\mid\sigma^k=\tau^{2^d}=1,
\tau\sigma\tau^{-1}=\sigma^{-1}\rangle$ where $d\geq 1, k\geq 3$, 
$n,k$: odd, and ${\rm gcd}\{n,k\}=1$.
\end{theorem}

When an extension $K/k$ is non-Galois and separable, 
we take the Galois closure $L/k$ of $K/k$ and 
put $G={\rm Gal}(L/k)$ and $H={\rm Gal}(L/K)$. 
Then we have:
\begin{theorem}[Endo {\cite[Theorem 2.1]{End11}}]
Assume that $G={\rm Gal}(L/k)$ is nilpotent. 
Then $R^{(1)}_{K/k}(\bG_m)$ is not 
retract $k$-rational.
\end{theorem}
\begin{theorem}[Endo {\cite[Theorem 3.1]{End11}}]\label{th15}
Assume that the Sylow subgroups of $G={\rm Gal}(L/k)$ are all cyclic. 
Then $R^{(1)}_{K/k}(\bG_m)$ is retract $k$-rational, 
and the following conditions are equivalent:\\
{\rm (i)} 
$R^{(1)}_{K/k}(\bG_m)$ is stably $k$-rational;\\
{\rm (ii)} 
$G=D_n$ with $n$ odd $(n\geq 3)$ 
or $G=C_m\times D_n$ where $m,n$ are odd, 
$m,n\geq 3$, ${\rm gcd}\{m,n\}=1$, and $H={\rm Gal}(L/K)\leq D_n$ is of order $2$;\\
{\rm (iii)} 
$H=C_2$ and $G\simeq C_r\rtimes H$, $r\geq 3$ odd, where 
$H$ acts non-trivially on $C_r$. 
\end{theorem}
\begin{theorem}[{Endo \cite[Theorem 4.1]{End11}, see also \cite[Remark 4.2]{End11}}]\label{thS}
\\Assume that ${\rm Gal}(L/k)=S_n$, 
$n\geq 3$, 
and ${\rm Gal}(L/K)=S_{n-1}$ is the stabilizer of one of the letters 
in $S_n$.\\
{\rm (1)}\ 
$R^{(1)}_{K/k}(\bG_m)$ is retract $k$-rational 
if and only if $n$ is a prime;\\
{\rm (2)}\ 
$R^{(1)}_{K/k}(\bG_m)$ is $($stably$)$ $k$-rational 
if and only if 
$n=3$.
\end{theorem}
\begin{theorem}[Endo {\cite[Theorem 4.4]{End11}}]\label{thA}
Assume that ${\rm Gal}(L/k)=A_n$, 
$n$ $\geq$ $4$, 
and ${\rm Gal}(L/K)=A_{n-1}$ is the stabilizer of one of the letters 
in $A_n$.\\
{\rm (1)}\ 
$R^{(1)}_{K/k}(\bG_m)$ is retract $k$-rational 
if and only if $n$ is a prime.\\
{\rm (2)}\ For some positive integer $t$, 
$[R^{(1)}_{K/k}(\bG_m)]^{(t)}$ is stably $k$-rational 
if and only if 
$n=5$, where 
$[R^{(1)}_{K/k}(\bG_m)]^{(t)}$ is the product of $t$ copies of 
$R^{(1)}_{K/k}(\bG_m)$. 
\end{theorem}

A birational classification of the algebraic $k$-tori of dimension $4$ 
is given by Hoshi and Yamasaki \cite{HY}. 
Note that there are $710$ $\bZ$-classes forming $227$ $\bQ$-classes in $GL_4(\bZ)$. 
\begin{theorem}[{\cite[Theorem 1.8]{HY}}]\label{th1}
Let $L/k$ be a Galois extension and $G\simeq 
{\rm Gal}(L/k)$ be a finite subgroup of $GL_4(\bZ)$ 
which acts on $L(x_1,x_2,x_3,x_4)$ via $(\ref{acts})$. \\
{\rm (i)} 
$L(x_1,x_2,x_3,x_4)^G$ is stably $k$-rational 
if and only if 
$G$ is conjugate to one of the $487$ groups which are not in 
{\rm \cite[Tables $2$, $3$ and $4$]{HY}}.\\
{\rm (ii)} 
$L(x_1,x_2,x_3,x_4)^G$ is not stably but retract $k$-rational 
if and only if $G$ is conjugate to one of the $7$ groups which are 
given as in {\rm \cite[Table $2$]{HY}}.\\
{\rm (iii)} 
$L(x_1,x_2,x_3,x_4)^G$ is not retract $k$-rational if and only if 
$G$ is conjugate to one of the $216$ groups which are given as 
in {\rm \cite[Tables $3$ and $4$]{HY}}.
\end{theorem}

Let $F_{20}$ be the Frobenius group of order $20$. 
By Theorem \ref{th1}, we have:
\begin{theorem}[{\cite[Theorem 1.9]{HY}}]\label{th11}
Let $K/k$ be a separable field extension of degree $5$ 
and $L/k$ be the Galois closure of $K/k$. 
Assume that $G={\rm Gal}(L/k)$ is a transitive subgroup of $S_5$ 
which acts on $L(x_1,x_2,x_3,x_4)$ via $(\ref{acts})$, 
and $H={\rm Gal}(L/K)$ is the stabilizer of one 
of the letters in $G$.\\
{\rm (1)} $R_{K/k}^{(1)}(\bG_m)$ is stably $k$-rational if and only if 
$G\simeq C_5$, $D_5$ or $A_5$;\\
{\rm (2)} $R_{K/k}^{(1)}(\bG_m)$ is not stably but retract $k$-rational 
if and only if $G\simeq F_{20}$ or $S_5$.
\end{theorem}
Theorem \ref{th11} is already known except for the case of $A_5$ 
(see Theorems \ref{th13}, \ref{th15}, \ref{thS} and \ref{thA}). 
Stable $k$-rationality of $R_{K/k}^{(1)}(\bG_m)$ for the case $A_5$ 
is asked by S. Endo in \cite[Remark 4.6]{End11}.
By Theorems \ref{thA} and \ref{th11}, we get: 
\begin{corollary}[{\cite[Corollary 1.10]{HY}}]\label{corA}
Let $K/k$ be a non-Galois separable field extension 
of degree $n$ and $L/k$ be the Galois closure of $K/k$. 
Assume that ${\rm Gal}(L/k)=A_n$, $n\geq 4$, 
and ${\rm Gal}(L/K)=A_{n-1}$ is the stabilizer of one of the letters 
in $A_n$. Then 
$R^{(1)}_{K/k}(\bG_m)$ is stably $k$-rational 
if and only if 
$n=5$.
\end{corollary}

There are $6079$ $\bZ$-classes forming $955$ $\bQ$-classes in $GL_5(\bZ)$. 
A birational classification of the algebraic $k$-tori 
of dimension $5$ is given as follows:
\begin{theorem}[{\cite[Theorem 1.11]{HY}}]\label{th2}
Let $L/k$ be a Galois extension and $G\simeq 
{\rm Gal}(L/k)$ be a finite subgroup of $GL_5(\bZ)$ 
which acts on $L(x_1,x_2,x_3,x_4,x_5)$ via $(\ref{acts})$.\\
{\rm (i)} 
$L(x_1,x_2,x_3,x_4,x_5)^G$ is stably $k$-rational if and only if 
$G$ is conjugate to one of the $3051$ groups which are not in 
{\rm \cite[Tables $11$, $12$, $13$, $14$ and $15$]{HY}}.\\
{\rm (ii)} 
$L(x_1,x_2,x_3,x_4,x_5)^G$ is not stably but retract $k$-rational 
if and only if $G$ is conjugate to one of the $25$ groups which are given as 
in {\rm \cite[Table $11$]{HY}}.\\
{\rm (iii)} 
$L(x_1,x_2,x_3,x_4,x_5)^G$ is not retract $k$-rational if and only if 
$G$ is conjugate to one of the $3003$ groups which are given as 
in {\rm \cite[Tables $12$, $13$, $14$ and $15$]{HY}}.
\end{theorem}
\begin{theorem}[{\cite[Theorem 1.13]{HY}}]\label{th22}
Let $K/k$ be a separable field extension 
of degree $6$ and $L/k$ be the Galois closure of $K/k$. 
Assume that $G={\rm Gal}(L/k)$ is a transitive subgroup of $S_6$ 
which acts on $L(x_1,x_2,x_3,x_4,x_5)$ via $(\ref{acts})$, 
and $H={\rm Gal}(L/K)$ is the stabilizer of one 
of the letters in $G$. Then 
$R_{K/k}^{(1)}(\bG_m)$ is stably $k$-rational if and only 
if $G\simeq C_6$, $S_3$ or $D_6$. 
Moreover, if $R_{K/k}^{(1)}(\bG_m)$ is not stably $k$-rational, 
then it is not retract $k$-rational. 
\end{theorem}

In Theorems \ref{th1}, \ref{th11}, \ref{th2} and \ref{th22}, 
when the field is stably $k$-rational, 
we do not know whether $L(M)^G$ is $k$-rational 
except for few cases (see \cite[Chapter 2]{Vos98}). 

\subsection{Strategy of the proof: the flabby class $[M]^{fl}$ of a $G$-lattice $M$}
\label{sse32}

Let $M$ be a $G$-lattice. 
For the rationality problem of $L(M)^G$ over $k$, 
the flabby class $[M]^{fl}$ of $M$ 
plays crucial role as follows (see Voskresenskii \cite[Section 4.6]{Vos98}, 
Lorenz \cite[Section 9.5]{Lor05}): 

\begin{theorem}[Endo and Miyata, Voskresenskii, Saltman]\label{thEM73}
Let $L/k$ be a finite Galois extension with Galois group $G={\rm Gal}(L/k)$ 
and $M, M^\prime$ be $G$-lattices.\\
{\rm (1)} {\rm (Endo and Miyata \cite[Theorem 1.6]{EM73})}\\ 
$[M]^{fl}=0$ if and only if $L(M)^G$ is stably $k$-rational.\\
{\rm (2)} {\rm (Voskresenskii \cite[Theorem 2]{Vos74})}\\ 
$[M]^{fl}=[M^\prime]^{fl}$ if and only if $L(M)^G$ and $L(M^\prime)^G$ 
are stably $k$-isomorphic.\\
{\rm (3)} {\rm (Saltman \cite[Theorem 3.14]{Sal84b})}\\ 
$[M]^{fl}$ is invertible if and only if $L(M)^G$ is 
retract $k$-rational.
\end{theorem}

Unfortunately, the flabby class $[M]^{fl}$ 
is useful to verify the stable $k$-rationality and the retract $k$-rationality 
of $L(M)^G$ but useless to verify the $k$-rationality. 

In order to give the definition of the flabby class 
$[M]^{fl}$ of $G$-lattice $M$ (Definition \ref{defF}), 
we prepare some terminology. 

\begin{definition}
Let $M$ be a $G$-lattice.\\ 
(1) $M$ is called a {\it permutation} $G$-lattice if $M$ has a $\bZ$-basis
permuted by $G$, i.e. $M\simeq \oplus_{1\leq i\leq m}\bZ[G/H_i]$ 
for some subgroups $H_1,\ldots,H_m$ of $G$.\\
(2) $M$ is called a {\it stably permutation} $G$-lattice if $M\oplus P\simeq P^\prime$ 
for some permutation $G$-lattices $P$ and $P^\prime$.\\ 
(3) $M$ is called {\it invertible} if it is a direct summand of a permutation $G$-lattice, 
i.e. $P\simeq M\oplus M^\prime$ for some permutation $G$-lattice 
$P$ and a $G$-lattice $M^\prime$.\\ 
(4) $M$ is called {\it coflabby} if $H^1(H,M)=0$
for any subgroup $H$ of $G$.\\ 
(5) $M$ is called {\it flabby} if $\widehat H^{-1}(H,M)=0$ 
for any subgroup $H$ of $G$ where $\widehat H$ is the Tate cohomology. 
\end{definition}
It is not difficult to verify the following implications:
\begin{center}
permutation\ \ $\Rightarrow$\ \
stably permutation\ \ $\Rightarrow$\ \ 
invertible\ \ $\Rightarrow$\ \ 
flabby and coflabby. 
\end{center}
\begin{definition}[{see \cite[Section 1]{EM74}, \cite[Section 4.7]{Vos98}}]
Let $\cC(G)$ be the category of all $G$-lattices. 
Let $\cS(G)$ be the full subcategory of $\cC(G)$ of all permutation $G$-lattices 
and $\cD(G)$ be the full subcategory of $\cC(G)$ of all invertible $G$-lattices. 
Let 
\begin{align*}
\cH^i(G)=\{M\in \cC(G)\mid \widehat H^i(H,M)=0\ {\rm for\ any}\ H\leq G\}\ (i=\pm 1)
\end{align*}
be the class of ``$\widehat H^i$-vanish'' $G$-lattices 
where $\widehat H^i$ is the Tate cohomology. 
Then one has the inclusions 
$\cS(G)\subset \cD(G)\subset \cH^i(G)\subset \cC(G)$ $(i=\pm 1)$. 
\end{definition}

\begin{definition}[The commutative monoid $\cT(G)=\cC(G)/\cS(G)$]
We say that two $G$-lattices $M_1$ and $M_2$ are {\it similar} 
if there exist permutation $G$-lattices $P_1$ and $P_2$ such that 
$M_1\oplus P_1\simeq M_2\oplus P_2$. 
We denote the set of similarity classes $\cC(G)/\cS(G)$ by $\cT(G)$ 
and the similarity class of $M$ by $[M]$. 
$\cT(G)$ becomes a commutative monoid 
with respect to the sum $[M_1]+[M_2]:=[M_1\oplus M_2]$ 
and the zero $0=[P]$ where $P\in \cS(G)$. 
\end{definition}

\begin{definition}\label{defF}
For a $G$-lattice $M$, there exists a short exact sequence of $G$-lattices
$0 \rightarrow M \rightarrow P \rightarrow F \rightarrow 0$
where $P$ is permutation and $F$ is flabby which is called a 
{\it flasque resolution} of $M$ (see Endo and Miyata \cite[Lemma 1.1]{EM74}, 
Colliot-Th\'el\`ene and Sansuc \cite[Lemma 3]{CTS77}). 
The similarity class $[F]\in \cT(G)$ of $F$ is determined uniquely and is called 
{\it the flabby class} of $M$. 
We denote the flabby class $[F]$ of $M$ by $[M]^{fl}$. 
We say that $[M]^{fl}$ is invertible if $[M]^{fl}=[E]$ for some 
invertible $G$-lattice $E$. 
\end{definition}
%



%
\begin{theorem}[{Colliot-Th\'{e}l\`{e}ne and Sansuc \cite[Corollaire 1]{CTS77}}]
\label{thCTS77}
Let $G$ be a finite group. 
The following conditions are equivalent:\\
{\rm (i)} $[J_G]^{fl}$ is coflabby;\\
{\rm (ii)} any Sylow subgroup of $G$ is cyclic or generalized quaternion $Q_{4n}$ 
of order $4n$ $(n\geq 2)$;\\
{\rm (iii)} any abelian subgroup of $G$ is cyclic;\\
{\rm (iv)} $H^3(H,\bZ)=0$ for any subgroup $H$ of $G$.
\end{theorem}
\begin{theorem}[{Endo and Miyata \cite[Theorem 2.1]{EM82}}]\label{thEM82}
Let $G$ be a finite group. 
The following conditions are equivalent:\\
{\rm (i)} $\cH^1(G)\cap \cH^{-1}(G)=\cD(G)$, i.e. 
any flabby and coflabby $G$-lattice is invertible;\\
{\rm (ii)} $[J_G\otimes_\bZ J_G]^{fl}=[[J_G]^{fl}]^{fl}$ is invertible;\\
{\rm (iii)} any p-Sylow subgroup of $G$ is cyclic for odd $p$ and 
cyclic or dihedral $($including Klein's four group$)$ for $p=2$.
\end{theorem}
Note that $H^1(H,[J_G]^{fl})\simeq H^3(H,\bZ)$ for any subgroup $H$ of $G$ 
(see \cite[Theorem 7]{Vos70} and \cite[Proposition 1]{CTS77}) 
and $[J_G]^{fl}=[J_G\otimes_\bZ J_G]$ (see \cite[Section 2]{EM82}). 

For $G$-lattice $M$, it is not difficult to see 
\begin{center}
permutation\ \ $\Rightarrow$\ \ 
stably permutation\ \ $\Rightarrow$\ \ 
invertible\ \ $\Rightarrow$\ \ flabby and coflabby\\
\hspace*{8mm}$\Downarrow$\hspace*{31mm}$\Downarrow$~~~~~~~~~~~~~\\
\hspace*{18mm}$[M]^{fl}=0$ in $\cT(G)$\ \ $\Rightarrow$
\ $[M]^{fl}$ is invertible in $\cT(G)$.
\end{center}

The above implications in each step cannot be reversed. 
Swan \cite{Swa60} gave an example of $Q_8$-lattice 
$M$ of rank $8$ which is not permutation but stably 
permutation: $M\oplus \bZ\simeq \bZ[Q_8]\oplus \bZ$. 
This also indicates that the direct sum cancellation fails. 
Colliot-Th\'el\`ene and Sansuc \cite[Remarque R1]{CTS77}, 
\cite[Remarque R4]{CTS77} 
gave examples of $S_3$-lattice $M$ of rank $4$ which 
is not permutation but stably permutation: 
$M\oplus \bZ\simeq \bZ[S_3/\langle\sigma\rangle]
\oplus\bZ[S_3/\langle\tau\rangle]$ 
where $S_3=\langle\sigma,\tau\rangle$ and 
of $F_{20}$-lattice $[J_{F_{20}/C_4}]^{fl}$ 
of the Chevalley module $J_{F_{20}/C_4}$ of rank $4$ 
which is  not stably permutation but invertible 
(see also Theorem \ref{thE0} and Theorem \ref{th1M} (ii), (iv) and (v)). 
By Theorem \ref{th13} (1), Theorem \ref{thEM73} (2) and Theorem \ref{thCTS77}, 
the flabby class $[J_{Q_8}]^{fl}$ of the Chevalley module $J_{Q_8}$ 
of rank $7$ is not invertible but flabby and coflabby 
(we may take $[J_{Q_8}]^{fl}$ of rank $9$, see \cite[Example 7.3]{HY}). 
The inverse direction of the vertical implication holds 
if $M$ is coflabby (see \cite[Lemma 2.11]{HY}).


Theorem \ref{th13} (resp. Theorem \ref{th15}, Theorem \ref{thS}, Theorem \ref{thA}) 
can be obtained by Theorem \ref{thEM7415} and Theorem \ref{thEM74M} 
(resp. Theorem \ref{thE0}, Theorem \ref{thSS}, Theorem \ref{thAA}) 
by using the interpretation as in Theorem \ref{thEM73}.

\begin{theorem}[{Endo and Miyata \cite[Theorem 1.5]{EM74}}]\label{thEM7415}
Let $G$ be a finite group. The following conditions are equivalent:\\
{\rm (i)} $[J_G]^{fl}$ is invertible;\\
{\rm (ii)} all the Sylow subgroups of $G$ are cyclic;\\
{\rm (iii)} $\cH^{-1}(G)=\cH^1(G)=\cD(G)$, i.e. any flabby $($resp. coflabby$)$ 
$G$-lattice is invertible.
\end{theorem}
\begin{theorem}[{
\cite[Theorem 2.3]{EM74}, see also 
\cite[Proposition 3]{CTS77}}]\label{thEM74M}
Let $G$ be a finite group. The following conditions are equivalent:\\
{\rm (i)} $[J_G]^{fl}=0$;\\
{\rm (ii)} $[J_G]^{fl}$ is of finite order in $\cT(G)$;\\
{\rm (iii)} all the Sylow subgroups of $G$ are cyclic and 
$H^4(G,\bZ)\simeq\widehat H^0(G,\bZ)$;\\
{\rm (iv)} $G=C_m$ or $G=C_n\times \langle\sigma,\tau\mid\sigma^k=\tau^{2^d}=1,
\tau\sigma\tau^{-1}=\sigma^{-1}\rangle$ where $d\geq 1, k\geq 3$, 
$n,k$: odd, and ${\rm gcd}\{n,k\}=1$;\\
{\rm (v)} $G=\langle s,t\mid s^m=t^{2^d}=1, tst^{-1}=s^r, m: odd,\ 
r^2\equiv 1\pmod{m}\rangle.$
\end{theorem}
\begin{theorem}[{Endo \cite[Theorem 3.1]{End11}}]\label{thE0}
Let $G$ be a non-abelian group. 
Assume that Sylow subgroups of $G$ are all cyclic.
Let $H$ be a non-normal subgroup of $G$ which contains no normal subgroup 
of $G$ except $\{1\}$. 
$($By Theorem \ref{thEM7415}, $[J_{G/H}]^{fl}$ is invertible.$)$ 
The following conditions are equivalent:\\
{\rm (i)} $[J_{G/H}]^{fl}=0$;\\
{\rm (ii)} $[J_{G/H}]^{fl}$ is of finite order in $\cT(G)$;\\
{\rm (iii)} 
$G=D_n$ with $n$ odd $(n\geq 3)$ 
or $G=C_m\times D_n$ where $m,n$ are odd, 
$m,n\geq 3$, ${\rm gcd}\{m,n\}=1$, and $H\leq D_n$ is of order $2$;\\
{\rm (iv)} 
$H=C_2$ and $G\simeq C_r\rtimes H$, $r\geq 3$ odd, where 
$H$ acts non-trivially on $C_r$. 
\end{theorem}

\begin{theorem}[{Endo \cite[Theorem 4.3]{End11}, see also \cite[Remark 4.2]{End11}}]\label{thSS}\\
Let $n\geq 3$ be an integer.\\
{\rm (1)} $[J_{S_n/S_{n-1}}]^{fl}$ is invertible if and only if $n$ is a prime.\\
{\rm (2)} $[J_{S_n/S_{n-1}}]^{fl}=0$ if and only if $n=3$.\\
{\rm (3)} $[J_{S_n/S_{n-1}}]^{fl}$ is of finite order in $\cT(G)$ 
if and only if $n=3$.
\end{theorem}
\begin{theorem}[{Endo \cite[Theorem 4.5]{End11}}]\label{thAA}
Let $n\geq 4$ be an integer.\\
{\rm (1)} $[J_{A_n/A_{n-1}}]^{fl}$ is invertible if and only if $n$ is a prime.\\
{\rm (2)} $[J_{A_n/A_{n-1}}]^{fl}$ is of finite order in $\cT(G)$ 
if and only if $n=5$.
\end{theorem}
%


We conclude this article, by introducing results on $G$-lattices of rank $4$ 
and $5$. 

\begin{definition}[The $G$-lattice $M_G$]\label{defMG} 
Let $G$ be a finite subgroup of $GL_n(\bZ)$. 
The $G$-lattice $M_G$ of rank $n$ 
is defined to be the $G$-lattice with a $\bZ$-basis $\{u_1,\ldots,u_n\}$ 
on which $G$ acts by $\sigma(u_i)=\sum_{j=1}^n a_{ij}u_j$ for any $
\sigma=[a_{ij}]\in G$. 
\end{definition}

For $2\leq n\leq 4$, the GAP code $(n,i,j,k)$ of a finite subgroup 
$G$ of $GL_n(\bZ)$ stands for 
the $k$-th $\bZ$-class of the $j$-th $\bQ$-class of 
the $i$-th crystal system of dimension $n$ as in 
\cite[Table 1]{BBNWZ78} and \cite{GAP}. 

A birational classification of the $k$-tori of dimension $4$ and $5$ 
(Theorem \ref{th1} and Theorem \ref{th2}) 
can be obtained by Theorem \ref{th1M} and Theorem \ref{th2M} respectively. 
\begin{theorem}[{\cite[Theorem 1.25]{HY}}]\label{th1M}
Let $G$ be a finite subgroup of $GL_4(\bZ)$ and 
$M_G$ be the $G$-lattice as in Definition \ref{defMG}.\\
{\rm (i)} 
$[M_G]^{fl}=0$ if and only if 
$G$ is conjugate to one of the $487$ groups which are not in 
{\rm \cite[{\rm Tables} $2$, $3$ and $4$]{HY}}.\\
{\rm (ii)} 
$[M_G]^{fl}$ is not zero but invertible 
if and only if $G$ is conjugate to one of the $7$ groups which are 
given as in {\rm \cite[{\rm Table} $2$]{HY}}.\\
{\rm (iii)} 
$[M_G]^{fl}$ is not invertible if and only if 
$G$ is conjugate to one of the $216$ groups which are given as 
in {\rm \cite[{\rm Tables} $3$ and $4$]{HY}}.\\
{\rm (iv)} 
$[M_G]^{fl}=0$ if and only if $[M_G]^{fl}$ is of finite order in $\cT(G)$.\\
{\rm (v)} For $G\simeq S_5$ of the GAP code $(4,31,5,2)$ in {\rm (ii)}, 
we have 
$
-[M_G]^{fl}=[J_{S_5/S_4}]^{fl}\neq 0.
$\\
{\rm (vi)} For $G\simeq F_{20}$ of the GAP code $(4,31,1,4)$ in {\rm (ii)}, 
we have 
$-[M_G]^{fl}=[J_{F_{20}/C_4}]^{fl}\neq 0$.
\end{theorem}

\begin{theorem}[{\cite[Theorem 1.26]{HY}}]\label{th2M}
Let $G$ be a finite subgroup of $GL_5(\bZ)$ and 
$M_G$ be the $G$-lattice as in Definition \ref{defMG}.\\
{\rm (i)} 
$[M_G]^{fl}=0$ if and only if 
$G$ is conjugate to one of the $3051$ groups which are not in 
{\rm \cite[{\rm Tables} $11$, $12$, $13$, $14$ and $15$]{HY}}.\\
{\rm (ii)} 
$[M_G]^{fl}$ is not zero but invertible 
if and only if $G$ is conjugate to one of the $25$ groups which are given as 
in {\rm \cite[{\rm Table} $11$]{HY}}.\\
{\rm (iii)} 
$[M_G]^{fl}$ is not invertible if and only if 
$G$ is conjugate to one of the $3003$ groups which are given as 
in {\rm \cite[{\rm Tables} $12$, $13$, $14$ and $15$]{HY}}.\\
{\rm (iv)} 
$[M_G]^{fl}=0$ if and only if $[M_G]^{fl}$ is of finite order in $\cT(G)$.
\end{theorem}


\begin{remark}
{\rm (1)} By the interpretation as in Theorem \ref{thEM73}, 
Theorem \ref{th1M} (v), (vi) claims that 
the corresponding two tori $T$ and $T'$ of dimension $4$ 
are not stably $k$-rational and are not stably $k$-isomorphic 
each other but the torus $T\times T'$ of dimension $8$ 
is stably $k$-rational.\\
{\rm (2)} When $[M]^{fl}$ is invertible, 
the inverse element of $[M]^{fl}$ is $-[M]^{fl}=[[M]^{fl}]^{fl}$. 
Hence Theorem \ref{th1M} (v) also claims that 
$[[M_G]^{fl}]^{fl}=[J_{S_5/S_4}]^{fl}$ and 
$[[J_{S_5/S_4}]^{fl}]^{fl}=[M_G]^{fl}$. 
\end{remark}

Finally, we give an application of the results in \cite{HY} 
which provides the smallest example exhibiting the failure of 
the Krull-Schmidt theorem for permutation $G$-lattices 
(see Dress's paper \cite[Proposition 9.6]{Dre73}): 
\begin{proposition}[{\cite[Proposition 6.7]{HY}}]
Let $D_6$ be the dihedral group of order $12$ and 
$\{1\}$, $C_2^{(1)}$, $C_2^{(2)}$, $C_2^{(3)}$, $C_3$, $V_4$, 
$C_6$, $S_3^{(1)}$, $S_3^{(2)}$ and $D_6$ be the 
conjugacy classes of subgroups of $D_6$. 
Then the following isomorphism of permutation $D_6$-lattices holds: 
\begin{align*}
& ~{} \bZ[D_6]\oplus\bZ[D_6/V_4]^{\oplus 2}\oplus\bZ[D_6/C_6]
\oplus\bZ[D_6/S_3^{(1)}]\oplus\bZ[D_6/S_3^{(2)}]\\
\simeq & ~{} \bZ[D_6/C_2^{(1)}]\oplus\bZ[D_6/C_2^{(2)}]
\oplus\bZ[D_6/C_2^{(3)}]\oplus\bZ[D_6/C_3]\oplus\bZ^{\oplus 2}.
\end{align*}
\end{proposition}


\begin{acknowledgment}
The author thanks the 
referee for very careful reading of the manuscript 
and for many helpful suggestions which greatly improved the clarity of 
this paper.
\end{acknowledgment}


\begin{thebibliography}{BBNWZ78}
%

\bibitem[Ber10]{Ber10} 
G. Berhuy, 
{\it An introduction to Galois cohomology and its applications}, 
With a foreword by Jean-Pierre Tignol. 
London Mathematical Society Lecture Note Series, 377. 
Cambridge University Press, Cambridge, 2010.

\bibitem[Bog88]{Bog88}
F. A. Bogomolov,
\textit{The Brauer group of quotient spaces by linear group actions},
Math. USSR Izv. \textbf{30} (1988) 455--485.

\bibitem[BB]{BB}
F. A. Bogomolov and C. B\"ohning, 
\textit{Isoclinism and stable cohomology of wreath products}, 
arXiv:1204.4747.

\bibitem[BBNWZ78]{BBNWZ78} 
H. Brown, R. B\"ulow, J. Neub\"user, 
H. Wondratschek and H. Zassenhaus. 
{\it Crystallographic Groups of Four-Dimensional Space}, 
John Wiley, New York, 1978. 

\bibitem[CHHK]{CHHK}
H. Chu, A. Hoshi, S.-J. Hu and M. Kang, 
\textit{Noether's problem for groups of order 243}, preprint.

\bibitem[CHKK10]{CHKK10}
H. Chu, S.-J. Hu, M. Kang and B. E. Kunyavskii, 
\textit{Noether's problem and the unramified Brauer groups 
for groups of order 64},
Int. Math. Res. Not. IMRN \textbf{2010}, 2329--2366.

\bibitem[CHKP08]{CHKP08}
H. Chu, S.-J. Hu, M. Kang and Y. G. Prokhorov,
\textit{Noether's problem for groups of order 32},
J. Algebra \textbf{320} (2008) 3022--3035.

\bibitem[CK01]{CK01} H. Chu and M. Kang, 
{\it Rationality of $p$-group actions}, 
J. Algebra \textbf{237} (2001) 673--690.

\bibitem[CTS77]{CTS77} J.-L. Colliot-Th\'{e}l\`{e}ne and J.-J. Sansuc, 
{\it La R-\'{e}quivalence sur les tores}, 
Ann. Sci. \'{E}cole Norm. Sup. (4) \textbf{10} (1977) 175--229. 

\bibitem[CTS07]{CTS07}
J.-L. Colliot-Th\'el\`ene and J.-J. Sansuc, 
\textit{The rationality problem for fields of invariants under linear
algebraic groups {\rm (}with special regards to the Brauer
groups{\rm )}}, in ``Proc. International Conference, Mumbai, 2004"
edited by V. Mehta, Narosa Publishing House, 2007.

\bibitem[Dra83]{Dra83}
P. K. Draxl, \textit{Skew fields}, London Math. Soc. Lecture Note
Series vol. 81, Cambridge Univ. Press, Cambridge, 1983.

\bibitem[Dre73]{Dre73} A. W. M. Dress, 
{\it Contributions to the theory of induced representations}, 
Algebraic K-theory, II: 
"Classical'' algebraic K-theory and connections with arithmetic 
(Proc. Conf., Battelle Memorial Inst., Seattle, Wash., 1972), pp.183--240. 
Lecture Notes in Math., Vol. 342, Springer, Berlin, 1973. 

\bibitem[End11]{End11} S. Endo, 
{\it The rationality problem for norm one tori}, 
Nagoya Math. J. \textbf{202} (2011) 83--106. 

\bibitem[EM73]{EM73} S. Endo and T. Miyata, 
{\it Invariants of finite abelian groups}, 
J. Math. Soc. Japan \textbf{25} (1973) 7--26. 

\bibitem[EM74]{EM74} S. Endo and T. Miyata, 
{\it On a classification of the function fields of algebraic tori}, 
Nagoya Math. J. \textbf{56} (1975) 85--104.

\bibitem[EM82]{EM82} S. Endo and T. Miyata, 
{\it Integral representations with trivial first cohomology groups}, 
Nagoya Math. J. \textbf{85} (1982) 231--240. 

\bibitem[Fis15]{Fis15} 
E. Fischer, 
{\it Die Isomorphie der Invariantenk\"orper der endlichen 
Abel'schen Gruppen linearer Transformationen}, 
Nachr. K\"onigl. Ges. Wiss. G\"ottingen (1915) 77--80.

\bibitem[GAP]{GAP} 
The GAP Group, GAP -- Groups, Algorithms, and Programming, 
Version 4.4.12; 2008. (http://www.gap-system.org).

\bibitem[GMS03]{GMS03}
S. Garibaldi, A. Merkurjev and J-P. Serre, 
\textit{Cohomological invariants in Galois cohomology}, 
AMS Univ. Lecture Series, vol.28, Amer. Math. Soc., Providence, RI, 2003.

\bibitem[GS06]{GS06} 
P. Gille and T. Szamuely, 
{\it Central simple algebras and Galois cohomology}, 
Cambridge Studies in Advanced Mathematics, 101. 
Cambridge University Press, Cambridge, 2006.

\bibitem[Haj87]{Haj87}
M. Hajja, 
\textit{Rationality of finite groups of monomial automorphisms of $K(x,y)$},
J. Algebra \textbf{109} (1987) 46--51.

\bibitem[HK92]{HK92}
M. Hajja and M. Kang,
\textit{Finite group actions on rational function fields}, 
J. Algebra \textbf{149} (1992) 139--154.

\bibitem[HK94]{HK94}
M. Hajja and M. Kang,
\textit{Three-dimensional purely monomial group actions},
J. Algebra \textbf{170} (1994) 805--860.

\bibitem[HS64]{HS64}
M. Hall Jr. and J. K. Senior, \textit{The groups of order $2^n$
{\rm (}$n\le 6${\rm )}}, Macmillan, New York, 1964.

\bibitem[HKKi]{HKKi} 
A. Hoshi, M. Kang and H. Kitayama, 
{\it Quasi-monomial actions and some 4-dimensional 
rationality problems}, arXiv:1201.1332.

\bibitem[HKKu]{HKKu} 
A. Hoshi, M. Kang and B. E. Kunyavskii, 
{\it Noether's problem and unramified Brauer groups}, 
to appear in Asian J. Math. Preprint version: arXiv:1202.5812.

\bibitem[HKY11]{HKY11}
A. Hoshi, H. Kitayama and A. Yamasaki,
\textit{Rationality problem of three-dimensional monomial group actions},
J. Algebra \textbf{341} (2011) 45--108.

\bibitem[HR08]{HR08}
A. Hoshi and Y. Rikuna,
\textit{Rationality problem of three-dimensional purely monomial 
group actions: the last case}, 
Math. Comp. \textbf{77} (2008) 1823--1829.

\bibitem[HY]{HY} 
A. Hoshi and A. Yamasaki, 
{\it Rationality problem for algebraic tori}, arXiv:1210.4525.

\bibitem[Jam80]{Jam80}
R. James, \textit{The groups of order $p^6$ {\rm (}$p$ an odd
prime{\rm )}}, Math. Comp. \textbf{34} (1980) 613--637.

\bibitem[JNOB90]{JNOB90}
R. James, M. F. Newman and E. A. O'Brien, \textit{The groups of
order 128}, J. Algebra \textbf{129} (1990) 136--158.

\bibitem[JLY02]{JLY02} 
C. U. Jensen, A. Ledet and N. Yui, 
{\it Generic polynomials}, 
Constructive aspects of the inverse Galois problem. 
Mathematical Sciences Research Institute Publications, 45. 
Cambridge University Press, Cambridge, 2002. 

\bibitem[Kan90]{Kan90}
M. Kang, \textit{Constructions of Brauer-Severi varieties and norm
hypersurfaces}, Canadian J. Math. \textbf{42} (1990) 230--238.

\bibitem[Kan04]{Kan04} 
M. Kang, 
{\it Rationality problem of $\rm GL\sb 4$ group actions}, 
Adv. Math. \textbf{181} (2004) 321--352.

\bibitem[Kan05]{Kan05} 
M. Kang, {\it Some group actions on $K(x\sb 1,x\sb 2,x\sb 3)$}, 
Israel J. Math. \textbf{146} (2005) 77--92.

\bibitem[Kan12]{Kan12}
M. Kang, \textit{Retract rational fields}, 
J. Algebra \textbf{349} (2012) 22--37.

\bibitem[Kar87]{Kar87}
G. Karpilovsky, \textit{The Schur Multiplier}, London Math. Soc.
Monographs, vol.2, Oxford Univ. Press, 1987.

\bibitem[KMRT98]{KMRT98}
M.-A. Knus, A. Merkurjev, M. Rost and J.-P. Tignol, 
{\it The book of involutions}, With a preface in French by J. Tits. 
American Mathematical Society Colloquium Publications, 44. 
American Mathematical Society, Providence, RI, 1998.

\bibitem[Kun54]{Kun54}
H. Kuniyoshi, 
{\it On purely-transcendency of a certain field}, 
Tohoku Math. J. (2) \textbf{6} (1954) 101--108.

\bibitem[Kun55]{Kun55}
H. Kuniyoshi, {\it On a problem of Chevalley}, 
Nagoya Math. J. \textbf{8} (1955) 65--67.

\bibitem[Kun56]{Kun56}
H. Kuniyoshi, {\it Certain subfields of rational function fields}, 
Proceedings of the international symposium on algebraic number theory, 
Tokyo \& Nikko, 1955, 241--243, Science Council of Japan, Tokyo, 1956. 

\bibitem[Kun90]{Kun90}
B. E. Kunyavskii, 
{\it Three-dimensional algebraic tori}, (Russian) 
Translated in Selecta Math. Soviet. \textbf{9} (1990) 1--21. 
Investigations in number theory (Russian), 90--111, 
Saratov. Gos. Univ., Saratov, 1987. 

\bibitem[Kun10]{Kun10}
B. E. Kunyavskii, 
{\it The Bogomolov multiplier of finite simple groups}, 
Cohomological and geometric approaches to rationality problems, 
209--217, Progr. Math., 282, 
Birkh\"auser Boston, Inc., Boston, MA, 2010.

\bibitem[Len74]{Len74} 
H. W. Lenstra, Jr., 
\textit{Rational functions invariant 
under a finite abelian group}, Invent. Math. \textbf{25} (1974) 299--325.

\bibitem[Lor05]{Lor05}
M. Lorenz, 
{\it Multiplicative invariant theory}, 
Encyclopaedia of Mathematical Sciences, 135. 
Invariant Theory and Algebraic Transformation Groups, VI. 
Springer-Verlag, Berlin, 2005.

\bibitem[MT86]{MT86}
Y. I. Manin and M. A. Tsfasman, 
\textit{Rational varieties: algebra, geometry and arithmetric}, 
Russian Math. Survey \textbf{41} (1986) 51--116.

\bibitem[Mas55]{Mas55} K. Masuda, {\it On a problem of Chevalley}, 
Nagoya Math. J. \textbf{8} (1955) 59--63.

\bibitem[Mas68]{Mas68} K. Masuda, {\it Application of theory of the 
group of classes of projective modules to existence problem of 
independent parameters of invariant}, 
J. Math. Soc. Japan \textbf{20} (1968) 223--232. 

\bibitem[Mor12]{Mor12}
P. Moravec, \textit{Unramified Brauer groups of finite and
infinite groups}, Amer. J. Math. \textbf{134} (2012) 1679--1704.

\bibitem[Mor]{Mor} 
P. Moravec, \textit{Unramified Brauer groups and isoclinism}, 
arXiv:1203.2422. 

\bibitem[Pop98]{Pop98} 
S. Yu. Popov, {\it Galois lattices and their birational invariants}. 
(Russian) Vestn. Samar. Gos. Univ. Mat. Mekh. Fiz. Khim. Biol. 1998, 
no. 4, 71--83.

\bibitem[Pro10]{Pro10}
Y. G. Prokhorov, \textit{Fields of invariants of finite linear
groups}, in ``Cohomological and geometric approaches to
rationality problems", edited by F. Bogomolov and Y. Tschinkel,
Progress in Math. vol. 282, Birkh\"auser, Boston, 2010.

\bibitem[Roq63]{Roq63}
P. Roquette, \textit{On the Galois cohomology of the projective
linear group and its applications to the construction of generic
splitting fields of algebras}, Math. Ann. \textbf{150} (1963) 411--439.

\bibitem[Roq64]{Roq64}
P. Roquette, \textit{Isomorphisms of generic splitting fields of
simple algebras}, J. Reine Angew. Math. \textbf{214/215} (1964) 207--226.

\bibitem[Sal82a]{Sal82a}
D. J. Saltman, \textit{Generic Galois extensions and problems in
field theory}, Adv. Math. \textbf{43} (1982) 250--283.

\bibitem[Sal82b]{Sal82b}
D. J. Saltman, 
{\it Generic structures and field theory}, 
Algebraists' homage: papers in ring theory and related topics (New Haven, Conn., 1981), 
pp. 127--134, Contemp. Math., 13, Amer. Math. Soc., Providence, R.I., 1982. 

\bibitem[Sal84a]{Sal84a}
D. J. Saltman,
\textit{Noether's problem over an algebraically closed field},
Invent. Math. \textbf{77} (1984) 71--84.

\bibitem[Sal84b]{Sal84b} 
D. J. Saltman, 
{\it Retract rational fields and cyclic Galois extensions}, 
Israel J. Math. \textbf{47} (1984) 165--215. 

\bibitem[Sal85]{Sal85}
D. J. Saltman, \textit{The Brauer group and the center of generic matrices}, 
J. Algebra \textbf{97} (1985) 53--67. 

\bibitem[Sal87]{Sal87}
D. J. Saltman, \textit{Multiplicative field invariants}, 
J. Algebra \textbf{106} (1987) 221--238.

\bibitem[Sal90a]{Sal90a}
D. J. Saltman, \textit{Twisted multiplicative field invariants,
Noether's problem and Galois extensions}, 
J. Algebra \textbf{131} (1990) 535--558.

\bibitem[Sal90b]{Sal90b}
D. J. Saltman,
\textit{Multiplicative field invariants and the Brauer group},
J. Algebra \textbf{133} (1990) 533--544.

\bibitem[Sal00]{Sal00} 
D. J. Saltman, 
{\it A nonrational field, answering a question of Hajja}, 
Algebra and number theory, Lecture Notes in Pure and Appl. Math., 
208, 263--271, Dekker, New York, 2000. 

\bibitem[Ser79]{Ser79}
J-P. Serre, \textit{Local fields}, Springer GTM vol. 67,
Springer-Verlag, Berlin, 1979.

\bibitem[Ser02]{Ser02} 
J-P. Serre, {\it Galois cohomology}, 
Translated from the French by Patrick Ion and revised by the author. 
Corrected reprint of the 1997 English edition. 
Springer Monographs in Mathematics. Springer-Verlag, Berlin, 2002

\bibitem[Swa60]{Swa60} 
R. G. Swan, 
{\it Induced Representations and Projective Modules}, 
Ann. Math. \textbf{71} (1960) 552--578.

\bibitem[Swa69]{Swa69} 
R. G. Swan, 
{\it Invariant rational functions and a problem of Steenrod}, 
Invent. Math. \textbf{7} (1969) 148--158.

\bibitem[Swa83]{Swa83}
R. G. Swan, \textit{Noether's problem in Galois theory}, in ``Emmy
Noether in Bryn Mawr", edited by B. Srinivasan and J. Sally,
Springer-Verlag, Berlin, 1983, pp. 21--40.

\bibitem[Vos67]{Vos67}
V. E. Voskresenskii, \textit{On two-dimensional algebraic tori II}
Math. USSR Izv. \textbf{1} (1967) 691--696.

\bibitem[Vos70]{Vos70} V. E. Voskresenskii, 
{\it Birational properties of linear algebraic groups}, (Russian) 
Izv. Akad. Nauk SSSR Ser. Mat. \textbf{34} (1970) 3--19; 
translation in Math. USSR-Izv. \textbf{4} (1970) 1--17.

\bibitem[Vos74]{Vos74} V. E. Voskresenskii, 
{\it Stable equivalence of algebraic tori}, (Russian) 
Izv. Akad. Nauk SSSR Ser. Mat. \textbf{38} (1974) 3--10; 
translation in Math. USSR-Izv. \textbf{8} (1974) 1--7.

\bibitem[Vos98]{Vos98} V. E. Voskresenskii, 
{\it Algebraic groups and their birational invariants}, 
Translated from the Russian manuscript by Boris Kunyavskii, 
Translations of Mathematical Monographs, 179. 
American Mathematical Society, Providence, RI, 1998.

\bibitem[Yam12]{Yam12}
A. Yamasaki, \textit{Negative solutions to three-dimensional
monomial Noether problem}, 
J. Algebra \textbf{370} (2012) 46--78.
\end{thebibliography}
\end{document}